\documentclass[11pt]{article}

\usepackage{array, amsmath, amsfonts, amsthm, amssymb, subcaption}
\usepackage{authblk, color, soul, bm, graphicx, empheq, enumitem, bbold}
\usepackage[dvipsnames]{xcolor}
\usepackage[margin = 1.5cm, font = small, labelfont = bf, labelsep = period]{caption}
\usepackage{xspace}
\usepackage[utf8]{inputenc}
\usepackage[T1]{fontenc}
\usepackage[round]{natbib}
\usepackage[ruled]{algorithm2e}
\usepackage{multirow}
\usepackage{mathrsfs}
\usepackage{comment}
\usepackage{eurosym}
\PassOptionsToPackage{hyphens}{url}\usepackage{hyperref}
\usepackage{accents}
\usepackage{booktabs}
\usepackage{mathtools}
\usepackage{multirow}
\usepackage{nicefrac}
\usepackage{cases}

\usepackage{fullpage}[2cm]

\usepackage{cleveref}

\usepackage[colorinlistoftodos]{todonotes}
\usepackage{comment}

\usepackage{pgfplots}

\definecolor{dodgerblue}{rgb}{0.118,0.565,1}

\theoremstyle{plain}
\newtheorem{Th}{Theorem}
\newtheorem{Prop}{Proposition}

\theoremstyle{definition}
\newtheorem{Ex}{Example}
\newtheorem{Rmk}{Remark}
\newtheorem{Ass}{Assumption}
\crefname{Ass}{Assumption}{Assumptions}
\crefname{Prop}{Proposition}{Propositions}
\crefname{lem}{Lemma}{Lemmas}

\newtheorem{lem}{Lemma}

\newcommand{\eg}{\textit{e.g.}}
\newcommand{\ie}{\textit{i.e.}}

\newcommand{\set}[1]{\mathcal{#1}}
\newcommand{\subj}{\rm{s.t.}}

\newcommand{\ubar}[1]{\underaccent{\bar}{#1}}

\newcommand{\D}{\mathcal{D}}

\newcommand{\T}{\mathcal{T}}

\newcommand{\R}{\mathbb R}

\DeclareMathOperator{\E}{\mathbb{E}}

\begin{document}
\begin{center}
\LARGE
    Frequency Regulation with Storage: On Losses and Profits \\
    \vspace{0.25cm}
\normalsize
Dirk Lauinger\textsuperscript{1}, Fran\c{c}ois Vuille\textsuperscript{2}, Daniel Kuhn\textsuperscript{3}

\footnotesize
    \textsuperscript{1}Massachusetts Institute of Technology, USA, \href{mailto:lauinger@mit.edu}{lauinger@mit.edu} 
    
    \textsuperscript{2}Etat de Vaud, Switzerland, \href{mailto:francois.vuille@vaud.ch}{francois.vuille@vaud.ch} 
    
    \textsuperscript{3}Ecole polytechnique fédérale de Lausanne, Switzerland, \href{mailto:daniel.kuhn@epfl.ch}{daniel.kuhn@epfl.ch}
\end{center}

\noindent Low-carbon societies will need to store vast amounts of electricity to balance intermittent generation from wind and solar energy, for example, through frequency regulation. Here, we derive an analytical solution to the decision-making problem of storage operators who sell frequency regulation power to grid operators and trade electricity on day-ahead markets. Mathematically, we treat future frequency deviation trajectories as functional uncertainties in a receding horizon robust optimization problem. We constrain the expected terminal state-of-charge to be equal to some target to allow storage operators to make good decisions not only for the present but also the future. Thanks to this constraint, the amount of electricity traded on day-ahead markets is an implicit function of the regulation power sold to grid operators. The implicit function quantifies the amount of power that needs to be purchased to cover the expected energy loss that results from providing frequency regulation. We show how the marginal cost associated with the expected energy loss decreases with roundtrip efficiency and increases with frequency deviation dispersion. We find that the profits from frequency regulation over the lifetime of energy-constrained storage devices are roughly inversely proportional to the length of time for which regulation power must be committed.

\section{Introduction}

In November 1896 when America's first large-scale power plant at Niagara Falls began transmitting electricity to the city of Buffalo about 20~miles away, the electricity came in the form of alternating current. The reason was that ``\emph{unlike direct current, alternating current \emph{[could]} travel}'' \cite[p.~48]{GAB16}. Any form of current must travel at a high voltage, otherwise transportation losses are unreasonably large. Before consumption the voltage must be reduced again to make it safe for use. Until the advent of modern power electronics in the 1950s, only the voltage of alternating current could be easily changed with pairs of metal coils, so-called transformers. Today, most electricity is still transported as alternating current. The frequency of the current indicates the mismatch between electricity demand and supply because it is linked to the rotational speed of turbines in powerplants. If more power is produced than consumed, the turbines speed up and the frequency raises. If more power is consumed than produced, the turbines slow down and the frequency falls.

The purpose of frequency regulation is to insure electricity grids against unforeseen second-to-second supply and demand mismatches. Traditionally, this insurance has been provided by centralized power plants, often fired by fossil fuels. As wind and solar power plants replace fossil-fuel-fired power plants, power generation becomes more weather-dependent, which may increase the demand for and decrease the supply of frequency regulation. Electricity storage could help to fill the gap. Lithium-ion batteries, in particular, are considered a promising source of frequency regulation, thanks to their fast dynamics and rapid cost decline~\citep{MZ21}. 

Several studies~\citep{XW21, YH22} have claimed that it will be profitable to invest in lithium-ion batteries for frequency regulation in the near future. Such studies focused on battery costs and frequency regulation prices but often did not explicitly model the uncertain demand for frequency regulation. In our previous work \citep{V2GRO}, we accounted for EU delivery guarantees and found that they significantly reduce the profits from frequency regulation. We found the profits to depend on charging and discharging efficiencies, and on the dispersion of frequency deviations, which could be measured by the mean absolute deviation or the standard deviation.

In this work, we focus on two research questions:
\begin{enumerate}
    \item How do charging and discharging losses influence frequency regulation profits?
    \item What regulatory changes would make frequency regulation with storage more profitable?
\end{enumerate}
To answer these questions, we derive an analytical solution to a simplified decision-making problem of a storage operator selling frequency regulation in the continental European electricity grid and trading electricity on day-ahead wholesale or retail markets. We distill the following insights:
\begin{enumerate}
    \item Although frequency deviations vanish on average, the average power flow entering a storage device is nonincreasing in the amount of regulation power offered. The lower the roundtrip efficiency of the storage device and the higher the dispersion of the frequency deviations, the higher the dissipative losses incurred by the provision of regulation power. In a numerical case-study, we find that the losses amount to between $0.6$\% and $4.3$\% of the regulation power offered to grid operators. To our best knowledge, we are the first to quantify the expected dissipative losses as a function of roundtrip efficiency and frequency deviation dispersion. If one assumed a deterministic frequency deviation signal equal to its average value of zero, then one would be led to believe that there are no such losses.
    \item For storage devices that are constrained by their storage capacity and their initial state-of-charge, rather than their charging and discharging capacities, the profits from frequency regualtion over the lifetime of the devices are inversely proportional to the length of time for which frequency regulation must be committed.
    \item EU regulators can make frequency regulation with storage more profitable if they $(i)$ make it easier for small-scale storage operators to access wholesale electricity markets, or $(ii)$ establish an intraday market for frequency regulation.
\end{enumerate}

While the analytical solution provides general insights in storage operation, it requires several simplifying assumptions and may thus be of limited practical value for the operation of any specific storage device. The insights can, however, be used as a base line against which storage operators can compare the results of more detailed operational models. We describe the simplifying assumptions below and explain how they can be addressed in operational models. In our methodological development, we will introduce technical assumptions, which are very likely to hold in reality.

First, we assume that charging and discharging losses are constant. In reality, the losses depend on the instantaneous state-of-charge of the storage device, the charging and discharging power, and temperature. \cite{TN19} formulated a dynamic program to account for decision-dependent losses and found them to reduce regulation profits by 10\% to 20\% for lithium-ion batteries.

Second, we assume that storage operators only trade on day-ahead markets and that they must commit to constant market bids throughout the entire day. In practice, they could participate on intraday markets for the wholesale of electricity. \cite{lohndorf2023value} and \cite{EK23} developed multi-stage stochastic programs for participation in day-ahead and intraday markets. They found higher expected profits on intraday than day-ahead markets. Löhndorf and Wozabal pointed out that these profits can only be fully realized by assets with a low price impact.

Third, we assume no market power. This is reasonable for small storage assets, such as stationary batteries, whereas large storage assets, such as pooled pumped-hydro power plants may exert some market power. \cite{EB23} formulated a Stackelberg game to analyze the market power of monopolistic storage operators under various EU market regulations. \cite{virasjoki2020utility} developed a bi-level optimization model to assess the market power of electricity producers. They found that electricity producers can reduce storage investments if they exert market power.

Finally, we do not consider any battery degradation. \cite{KU18} have shown that battery degradation can be negligible for electric vehicles providing frequency regulation depending on the operating conditions. Other storage technologies may experience different degradation dynamics. \cite{BV19} formulated a stochastic mixed-integer linear program that accounts for battery degradation. They found that storage operators will want to protect themselves against deep battery discharges, which are particularly harmful to battery longevity, by limiting the amount of regulation power they sell to grid operators. 

In terms of methodological development, we use a constraint on the expected terminal state-of-charge, rather than a value function as in~\cite{V2GRO}, to steer the storage operator toward decisions that work well for both the present and the future. Expected value constraints have previously been investigated by \cite{DD03} to establish stochastic dominance of the second order, by \cite{PR11} in hypothesis testing, and by \cite{DF21} as a generalization of mean-variance models. Here, we will use them for analytical tractability. We show in Section~\ref{sec:det_pb} that the amount of power bought on electricity markets is an implicit function of the amount of regulation power sold to grid operators under the expected value constraint. This leads to a one-dimensional decision problem, which can be solved highly efficiently by bisection for general frequency deviation distributions (see Section~\ref{sec:cs}) and analytically for two- and three-point distributions (see Section~\ref{sec:as}).

To ease readability, we refer to generic electricity storage devices as batteries and relegate all essential proofs to Appendix~\ref{sec:ch3_apx_proofs} and all other proofs to the supplementary material (SM), which also contains additional description and analysis of the case studies presented in Section~\ref{sec:appl} and an extension of our model to elastic prices. Throughout the manuscript, we use the personal pronoun \emph{``she''} to refer to battery operators and the personal pronoun \emph{``he''} to refer to grid operators.

\paragraph*{Notation.}  We designate all random variables by tilde signs. Their realizations are denoted by the same symbols without tilde signs. For any $z \in \R$, we set $[z]^+ = \max\{z,0\}$ and $[z]^- = \max\{-z,0\}$. For any closed intervals $\set{T}, \set{U} \subseteq \mathbb{R}$, we define $\set{L}(\set{T}, \set{U})$ as the space of all Riemann integrable functions $f:\set{T} \to \set{U}$, and we denote the intersection of a set~$\set{B} \subseteq \set{L}(\T, \mathbb{R}) $ with $\set{L}(\T, \mathbb{R}_+)$ as $\set{B}^+$. For any signed function $\delta \in \set{L}(\T, \mathbb{R})$, we denote by~$\vert \delta \vert$ the absolute value function with $\vert \delta \vert (t) = \vert \delta(t) \vert$ for every~$t \in \T$. We use~$g'_-$ and~$g'_+$ to denote the left and right derivatives of a univariate proper convex function~$g:\mathbb{R} \to (-\infty, +\infty]$, respectively.
\section{Problem Description}
We study the decision problem of a battery operator who can sell frequency regulation power~$x^r \in \mathbb{R}_+$ to a grid operator and buy electric power~$x^b \in \mathbb{R}$ on a wholesale or retail market. We allow $x^b$ to be negative, in which case the amount~$\vert x^b \vert$ of power is sold. Both $x^r$ and $x^b$ are chosen ex ante and kept constant over a prescribed planning horizon $\T = [0, T]$ (\eg, the next day). 
At any time~$t \in \T$, the battery operator measures the normalized deviation~$\tilde \delta(t) \in [-1,1]$ of the uncertain instantaneous grid frequency~$\tilde \nu(t)$ from the nominal frequency~$\nu_0$ and must consume the amount~$x^b + \tilde \delta(t) x^r$ of power from the grid. Mathematically, the normalized frequency deviation at time~$t$ is given by the clipped ramp function
\begin{equation*}
    \tilde \delta(t) = \begin{cases}
    +1 & \text{if } \tilde \nu(t) > \nu_0 + \Delta \nu, \\
    \frac{\tilde \nu(t) - \nu_0}{\Delta \nu} & \text{if } \nu_0 - \Delta \nu \leq \tilde \nu(t) \leq \nu_0 + \Delta \nu, \\
    -1 & \text{if } \tilde \nu(t) < \nu_0 - \Delta \nu,
    \end{cases}
\end{equation*}
where~$\Delta \nu$ is the maximum frequency deviation against which the grid operator seeks protection. 

\begin{Rmk}
    Battery operators are required to measure the frequency deviation in near real-time, \eg, at least every 10~seconds in the European Union~\citep[art.~154(10)]{EU17}. Our continuous-time formulation is independent of any particular time discretization.\hfill $\Box$
\end{Rmk}

The remuneration for offering frequency regulation is twofold:  the power~$x^r$ set aside for frequency regulation is compensated at the availability price $\tilde p^a(t)$, and the regulation power $\tilde \delta(t) x^r$ actually delivered at time~$t$ is compensated at the delivery price~$\tilde p^d(t)$. 
The power~$x^b$ bought on the market is priced at $\tilde p^b(t)$. The expected cost over the planning horizon~$\T$ thus amounts to
\begin{equation*}
	\E \int_\T \tilde p^b(t) x^b - \left( \tilde p^a(t) - \tilde \delta(t) \tilde p^d(t) \right)x^r \, \mathrm{d}t.
\end{equation*}

The net power flow leaving the grid at time~$t \in \T$ is given by $x^b + \delta(t) x^r$. In the following, we find it useful to distinguish the charging power $y^+(x^b, x^r, \delta(t)) = [x^b + \delta(t) x^r]^+$ from the discharging power $y^-(x^b, x^r, \delta(t)) = [x^b + \delta(t) x^r]^-$. We assume that the charging power is bounded above by the charging capacity~$\bar y^+ \in \mathbb{R}_+$, and the discharging power is bounded above by the discharging capacity~$\bar y^- \in \mathbb{R}_+$. When the battery is charging ($y^+ > 0$), then only a fraction $\eta^+$ of the charging power enters the battery, where~$\eta^+ \in (0,1]$ represents the charging efficiency. The rest is dissipated during the charging process. Conversely, when the battery is discharging ($y^- > 0$), then a multiple~$\frac{1}{\eta^-}$ of the discharging power leaves the battery, where $\eta^- \in (0,1]$ represents the discharging efficiency. The battery state-of-charge at any time $t\in\T$ can thus be expressed as
\begin{equation}\label{eq:soc}
	y\big(x^b, x^r, \delta, y_0, t\big) = y_0 + \int_0^t \eta^+ y^+\big(x^b, x^r, \delta(t')\big) - \frac{1}{\eta^-} y^-\big(x^b, x^r, \delta(t')\big) \, \mathrm{d}t',
\end{equation}
where $y_0$ denotes the initial state-of-charge, and~$\delta \in \set{L}(\set{T}, [-1,1])$ is a given frequency deviation trajectory. Throughout the planning horizon, the battery state-of-charge must remain between~$0$ and the battery capacity~$\bar y > 0$. We assume from now on that $0 \leq y_0 \leq \bar y$. The following proposition establishes fundamental qualitative properties of the state-of-charge function~$y$.

\begin{Prop}\label{Prop:y}
	All else being equal, the battery state-of-charge $y(x^b,x^r,\delta,y_0,t)$ is concave and strictly increasing in $x^b$, concave in $x^r$, concave nondecreasing in $\delta$, and affine nondecreasing in $y_0$.
\end{Prop}
Proposition~\ref{Prop:y} slightly strengthens Proposition~1 by \cite{V2GRO}.

The battery may be used beyond the immediate planning horizon~$\T$ for selling more regulation power, for exchanging power on other electricity markets, or for supplying power to electric devices. The extent to which this is possible depends on the state-of-charge at the end of the immediate planning horizon. The value of any particular terminal state-of-charge~$y(x^b, x^r, \delta, y_0, T)$ could be captured by a reward-to-go function as in dynamic programming~\citep{bertsekas1976dp}. This allows the battery operator to trade off present and future costs when selecting~$x^b$ and~$x^r$. Another approach is to constrain the terminal state-of-charge to be close to some target~$y^\star$ that will guarantee satisfactory future performance. Both terminal costs and terminal constraints are widely studied in model predictive control~\citep{DM00}. Here, we will use a terminal constraint because it will allow us express~$x^b$ as an implicit function of~$x^r$, which facilitates the analytical solution.

The terminal state-of-charge is uncertain at time~$0$ when the battery operator selects~$x^b$ and~$x^r$ because it depends on the frequency deviation trajectory~$\delta$ during the planning horizon~$\T$. In fact, the battery operator can only be sure to meet a fixed target~$y^\star$ if she sells \emph{no} regulation power ($x^r = 0$), which shields her from the uncertainty of the frequency deviations. If she sells regulation power ($x^r > 0$), however, then all she can hope for is to reach a terminal state-of-charge that is close to the target~$y^\star$ on average. In the following, we will thus require that the terminal state-of-charge be equal to~$y^\star$ in expectation.
We emphasize that this constraint is \emph{not} dictated by physics but is simply a means to contain future operating costs, which are not modeled explicitly.

Throughout the planning horizon, the battery operator must be able to honor all market commitments for all reasonably likely frequency deviation trajectories~$\delta$. 
\cite{MA20} show that extreme frequency deviation trajectories are very uncommon. It would thus appear overly conservative to impose the charging, discharging, and battery state-of-charge constraints robustly for all possible frequency deviation trajectories. 
Inspired by applicable regulations by the \cite{EU17}, we assume instead that the battery operator must satisfy the constraints only for the frequency deviation trajectories in the uncertainty set
\begin{equation*}
	\D = \left\{
	\delta \in \mathcal{L}(\T, [-1, 1]): \int_{\T} \left \vert \delta(t) \right \vert \, \mathrm{d}t \leq \gamma
	\right\}
\end{equation*}
parametrized by the uncertainty budget~$\gamma \in (0,T]$. Note that~$\gamma$ represents the maximum amount of time for which a scenario~$\delta \in \D$ may adopt an extreme value~$\delta(t) \in \{-1,1\}$. Note also that~$\D$ can be seen as an extension of the budget uncertainty sets introduced by~\cite{DB04} to functional uncertainties. The following lemma establishes symmetry properties of~$\D$ that will allow us to reduce the decision problem of the battery operator to a deterministic optimization problem.
\begin{lem}\label{lem:D}
We have $\delta \in \D$ if and only if~$\vert \delta \vert \in \D^+$.
\end{lem}

In summary, the battery operator's decision problem is to select $x^b$ and $x^r$ so as to minimize expected costs while meeting the battery state-of-charge target~$y^\star$ in expectation and ensuring that the charger, discharger, and battery capacities are respected at all times and under all frequency deviation trajectories~$\delta \in \D$. This gives rise to the following optimization problem.

\begin{equation}
	\tag{R}
	\label{pb:P}
	\arraycolsep=4pt\def\arraystretch{1.2}
	\begin{array}{>{\displaystyle}c*3{>{\displaystyle}l}}
		\min_{x^b \in \R, \, x^r \in \R_+} & \multicolumn{3}{>{\displaystyle}l}{
		\E \int_\T \tilde p^b(t) x^b - \left( \tilde p^a(t) - \tilde \delta(t) \tilde p^d(t) \right)x^r \, \mathrm{d}t
		} \\
		\subj & y^+(x^b, x^r ,\delta(t)) &\leq \bar{y}^+ & \forall \delta \in \D,~ \forall t \in \T \\
		& y^-(x^b, x^r, \delta(t)) &\leq \bar{y}^- & \forall \delta \in \D,~ \forall t \in \T \\
		& y(x^b,x^r,\delta,y_0,t) &\leq \bar{y} & \forall \delta \in \D,~\forall t \in \T \\
		& y(x^b,x^r,\delta,y_0,t) &\geq 0 & \forall \delta \in \D,~ \forall t \in \T \\
		& \E \left[ y(x^b, x^r, \tilde \delta, y_0, T) \right] & = y^\star 
	\end{array}
\end{equation}
The battery operator only needs to insure frequency deviation trajectories in~$\D$. For trajectories outside of~$\D$, the battery operator has to deliver regulation power up to the smallest time instant~$t_\gamma$ such that~$\int_0^{t_\gamma} \vert \delta(t) \vert \, \text{d}t = \gamma$. At all time instants $t > t_\gamma$, the battery operator does not need to deliver any regulation power and may assume that~$\delta(t) = 0$. We thus assume that $\mathbb{P}[\tilde \delta \in \D] = 1$ throughout the rest of the paper. For later use, we note that any $\delta \in \D$ has a mean absolute deviation~$\frac{1}{T}\int_\T \vert \delta(t) \vert \, \text{d}t$ no greater than~$\frac{\gamma}{T}$.

For a fixed frequency deviation trajectory~$\delta$, the textbook approach to solving the deterministic counterpart of problem~\eqref{pb:P} is to first discretize the planning horizon into~$N$ periods and then introduce~$N$ binary variables expressing whether the battery is charging or discharging during the respective periods~\citep[p.~85]{JAT15}. This results in a large-scale mixed-integer linear program. In the remainder, we will show that the robustly-constrained optimization problem~\eqref{pb:P} is much easier to solve than its deterministic counterpart. In Section~\ref{sec:det_pb} we will first show that problem~\eqref{pb:P} is equivalent to a one-dimensional deterministic optimization problem. In Section~\ref{sec:cs}, we will then see that for realistic values of the roundtrip efficiency~$\eta^+\eta^-$ the search space can be reduced to merely three candidate solutions, which can all be computed highly efficiently by bisection. In Section~\ref{sec:as}, we will show that for specific distributions of the frequency deviations problem~\eqref{pb:P} can even be solved in closed form. These results imply that robustification reduces complexity and confirm our previous findings~\citep[Remark~4]{V2GRO}.
\section{Reduction to a Deterministic Optimization Problem}\label{sec:det_pb}

In order to simplify problem~\eqref{pb:P}, we first rewrite its objective function as an explicit linear function of the decision variables. Next, we show that the robust constraints are equivalent to deterministic linear constraints. Finally, we exploit the terminal state-of-charge constraint to express~$x^b$ as an implicit function of~$x^r$, which allows us to reformulate problem~\eqref{pb:P} only in terms~of~$x^r$.

Note first that the objective function of problem~\eqref{pb:P} can be expressed as~$T(c^b x^b - c^r x^r)$, where $c^b = \E \frac{1}{T} \int_\T \tilde p^b(t) \, \mathrm{d}t$ denotes the expected average market price of electricity, and $c^r = \E \frac{1}{T} \int_\T \tilde p^a(t) - \tilde \delta(t)\tilde p^d(t) \, \mathrm{d}t$ denotes the expected average price of regulation power. In the following, we will assume without much loss of generality that~$c^b > 0$ and~$c^r > 0$. 

We now show that the robust constraints are equivalent to deterministic linear constraints. This may be surprising because the state-of-charge is concave in the decision variables, implying that the upper bounds on the state-of-charge represent nonconvex constraints. Similarly, as the state-of-charge is concave in~$\delta$, finding the worst-case frequency deviation trajectories for the lower bounds on the state-of-charge amounts to solving a nonconvex optimization problem. In addition to these complications, the bounds on the state-of-charge also need to hold for all time instants in the planning horizon. 
In general, constraints with such properties are severely intractable.
Given that~$x^b$ and~$x^r$ do not depend on time, it is tempting to think that~$\delta$ can be restricted to a constant function of time without loss of generality. This restriction of the uncertainty set, however, relaxes the feasible set, and one can show that the relaxed feasible set contains decisions that are infeasible in practice. Indeed, averaging the real frequency deviation signals \emph{under}estimates the maximum state-of-charge and \emph{over}estimates the minimum state-of-charge~\citep[Example~1]{V2GRO}. 

Although the robust constraints of problem~\eqref{pb:P} seem intractable, we can reformulate them as deterministic linear constraints. This is possible because the worst-case frequency deviation trajectories and the worst-case time instants can be evaluated \emph{a priori}. The following proposition summarizes our results.

\begin{Prop}[Constraint reduction]\label{prop:cr}
	If~$0 \leq y_0 \leq \bar y$, then the following equivalences hold.
	\begin{equation*}
		\arraycolsep=4pt\def\arraystretch{1.2}
		\begin{array}{*6{>{\displaystyle}l}}
			(i) & y^+(x^b, x^r ,\delta(t)) &\leq \bar{y}^+ & \forall \delta \in \D,~ \forall t \in \T
			& \iff
			x^r + x^b \leq \bar y^+ 
			\\ (ii) &
			y^-(x^b, x^r, \delta(t)) &\leq \bar{y}^- & \forall \delta \in \D,~ \forall t \in \T 
			& \iff
			x^r - x^b \leq \bar y^-
			\\ (iii) &
			y(x^b,x^r,\delta,y_0,t) & \leq \bar{y} & \forall \delta \in \D,~\forall t \in \T
			& \iff
			x^r + \max\left\{ \frac{T}{\gamma} x^b, x^b \right\}
			& \leq \frac{\bar y - y_0}{\eta^+ \gamma}
			\\ (iv) &
			y(x^b,x^r,\delta,y_0,t) & \geq 0 & \forall \delta \in \D,~ \forall t \in \T
			& \iff
			x^r - \min\left\{ \frac{T}{\gamma} x^b, x^b \right\} 
			& \leq \frac{\eta^- y_0}{\gamma}
		\end{array}
	\end{equation*}
\end{Prop}

Proposition~\ref{prop:cr} is inspired by Theorems~1 and~2 by \cite{V2GRO}. The proof critically exploits the monotonicity properties of~$y$ established in Proposition~\ref{Prop:y} and the symmetry of the uncertainty set~$\D$ established in Lemma~\ref{lem:D}. The proof reveals that the upper bound on the charging power and the upper bound on the state-of-charge are valid for all frequency deviations signals $\delta \in \D$ and all time instants~$t \in \T$ if and only if they are valid for the time instants~$\gamma$ and~$T$ and for the particular frequency deviation signal~$\delta^{(+)}$, defined as $\delta^{(+)}(t) = 1$ if $t \leq \gamma$ and $\delta^{(+)}(t) = 0$ otherwise. Similarly, the upper bound on the discharging power and the lower bound on the state-of-charge are valid for all frequency deviations signals $\delta \in \D$ and all time instants~$t \in \T$ if and only if they are valid for the time instants~$\gamma$ and~$T$ and for the particular frequency deviation signal~$\delta^{(-)} = - \delta^{(+)}$. 

Intuitively, if~$x^b \geq 0$, then the maximum state-of-charge is achieved at time~$T$ by any nonnegative frequency deviation trajectory that exhausts the uncertainty budget, \ie, that has a cumulative deviation of~$\gamma$, such as~$\delta^{(+)}$. If~$x^b < 0$, then the maximum state-of-charge is achieved at time~$\gamma$ by the frequency deviation trajectory that exhausts the uncertainty budget as quickly as possible, \ie, that achieves a cumulative deviation of~$\gamma$ as soon as possible. Since~$\delta \in \set{L}(\set{T},[-1,1])$, the nonnegative signal that exhausts the uncertainty budget as quickly as possible is~$\delta^{(+)}$. As $\delta^{(+)}(\gamma) = 1$ and $\delta(t) \leq 1$ for all~$\delta \in \D$ and all~$t \in \T$, the maximum charging power is achieved at time~$\gamma$ by the frequency deviation trajectory~$\delta^{(+)}$. The intuition for the lower bound on the state-of-charge and the upper bound on the discharging power is similar.

In the following, we will exploit the constraint on the expected terminal state-of-charge to express~$x^b$ as an implicit function of~$x^r$. To this end, we first compress the stochastic process $\{ \tilde \delta(t) \}_{t \in \T}$ to a single random variable~$\tilde \xi = \tilde \delta(\tilde t)$, where~$\tilde t$ is a random time independent of $\tilde \delta$ that follows the uniform distribution on the planning horizon~$\set{T}$. The marginal probability distribution~$\mathbb{P}_\xi$ of~$\tilde \xi$ is defined through~$\mathbb{P}_\xi[\set{B}] = \mathbb{P}[\tilde \xi \in \set{B}]$ for every Borel set~$\set{B} \subseteq \mathbb{R}$. One readily verifies that
\begin{equation*}
    \mathbb{P}_\xi[\set{B}]
    = \mathbb{P}[\tilde \xi \in \set{B}]
    = \mathbb{E} \left[ \mathbb{P} \left[\tilde \delta(\tilde t) \in \set{B} \, \big \vert \, \tilde t \right] \right]
    = \frac{1}{T} \int_\T \mathbb{P} \left[ \tilde \delta(\tilde t) \in \set{B} \, \big \vert \, \tilde t = t \right] \, \mathrm{d}t
    = \frac{1}{T} \int_\T \mathbb{P}\left[ \tilde \delta(t) \in \set{B} \right] \, \mathrm{d}t
\end{equation*}
for every Borel set~$\set{B} \subseteq \set{T}$, where the third and the fourth equalities hold because~$\tilde t$ is uniformly distributed on~$\set{T}$ and because~$\tilde t$ is independent of~$\tilde \delta(t)$ for every~$t \in \set{T}$, respectively. Historical frequency deviation data suggest that the marginal distribution of~$\tilde \xi$ is symmetric around zero (see Figure~\ref{fig:F_sig} in SM \ref{sm:delta}).
From now on, we will thus make the following assumption.
\begin{Ass}[Symmetry]\label{ass:phi}
	We have $\mathbb{P}_{\xi}[\set{B}] = \mathbb{P}_{\xi}[-\set{B}]$ for every Borel set~$\set{B} \subseteq \set{T}$.
\end{Ass}

If the stochastic process~$\{ \tilde \delta(t) \}_{t \in \T}$ is stationary, then the marginal distribution of~$\tilde \delta(t)$ coincides with~$\mathbb{P}_\xi$ for every~$t \in \T$. Based on data from the UK and the continental European electricity grid, \cite{MA20} find that frequency deviations are stationary on timescales longer than 24~hours. The 24~hour threshold coincides with the typical length of planning horizons for frequency regulation. This suggests that~$\mathbb{P}_\xi$ does not change from one planning horizon to the next and can thus be estimated from historical data.

In the following, we define $F:\R \to [0,1]$ as the cumulative distribution function corresponding to~$\mathbb{P}_\xi$, and we define~$\varphi:\R \to \R_+$ as the antiderivative of~$F$ with~$\varphi(-1) = 0$. For short, we will refer to~$\varphi$ as the super-cumulative distribution function. We are now ready to investigate the expected terminal state-of-charge as a function of~$x^b$ and~$x^r$. To this end, we define $\eta_d = \frac{1}{\eta^-} - \eta^+$.

\begin{Prop}[Properties of the expected terminal state-of-charge]\label{prop:esoc}
	The expected terminal state-of-charge is continuous, 	jointly concave in~$x^b$ and~$x^r$, strictly increasing and unbounded above in~$x^b$, and nonincreasing in~$x^r$. In particular, it is given by 
	\begin{equation}\label{eq:esoc}
	\E \left[ y(x^b, x^r, \tilde \delta, y_0, T) \right] = 
	    y_0 + T \left( \eta^+ x^b - \eta_d x^r \varphi\left( - \frac{x^b}{x^r}\right) \right)
	    ~~\forall (x^b, x^r) \in \mathbb{R} \times \mathbb{R}_{+}.
	\end{equation}
\end{Prop}
Note that~$x^r \varphi(-\frac{x^b}{x^r})$ represents the perspective function of~$\varphi(-x^b)$, which is jointly convex in~$x^b$ and~$x^r$ because~$\varphi$ is convex~\citep[p.~89]{SB04}. For~$x^r = 0$, the perspective function is defined as~$\lim_{x^r \to 0^+} x^r \varphi(-\frac{x^b}{x^r})$ and coincides thus with $[x^b]^-$. 

\begin{figure}[t!]
	\centering
	\includegraphics[width=1\linewidth]{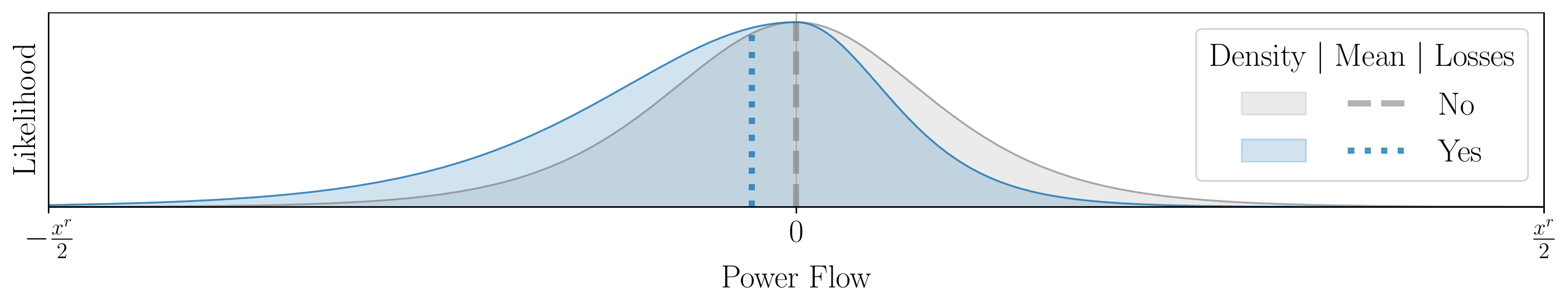}
	\caption{Distribution of the power flow $(\eta^+[\delta(t)]^+ - \frac{1}{\eta^-}[\delta(t)]^-)x^r$ entering the battery at any time~$t$ for $x^b = 0$.}
	\label{fig:losses_dist}
\end{figure}

We know from Proposition~\ref{Prop:y} that the battery state-of-charge is strictly increasing in~$x^b$, and thus it is unsurprising that its expected value is also strictly increasing in~$x^b$. We emphasize, however, that the state-of-charge may display a complicated nonmonotonic dependence on~$x^r$. Nevertheless, Proposition~\ref{prop:esoc} reveals that the \emph{expected} state-of-charge is nonincreasing in~$x^r$, which means that, on average, providing frequency regulation causes energy losses and thereby discharges the battery. Even though the average frequency deviations vanish by virtue of Assumption~\ref{ass:phi}, frequency regulation fails to be energy-neutral unless $\eta^+ = 1$ and $\eta^- = 1$. We will explain this phenomenon by reasoning about the power flow entering the battery as opposed to the power flow exiting the electricity grid. If there are no losses, then the power flow entering the battery and the power flow exiting the electricity grid coincide with~$x^b + \delta(t) x^r$. They thus follow the same probability distribution as the frequency deviations with the mean value shifted from~0 to~$x^b$ and the standard deviation scaled by~$x^r$. Providing frequency regulation hence only increases the dispersion of the power flow entering the battery but does not affect its mean. In the general case, when $\eta^+, \eta^- < 1$, the power flow exiting the electricity grid follows the same probability distribution as before, but the power flow entering the battery follows a different probability distribution. In fact, charging losses \emph{compress} the positive part of the original distribution, while discharging losses \emph{stretch} the negative part of the original distribution. The losses thereby \emph{decrease} the average power flow entering the battery. The higher the dispersion of the power flow, the more pronounced the decrease. Since the dispersion increases in~$x^r$, the average power flow entering the battery and, by extension, the expected terminal state-of-charge of the battery decrease in~$x^r$. Figure~\ref{fig:losses_dist} visualizes the distribution of the power flow entering the battery for~$x^b = 0$ with and without losses.

The monotonicity properties of the expected terminal state-of-charge established in Proposition~\ref{prop:esoc} imply that the last constraint of problem~\eqref{pb:P} determines~$x^b$ as an implicit function of~$x^r$. Instead of reasoning about the state-of-charge of the battery directly, we will continue to reason about the power flow entering the battery, which is independent of the initial state-of-charge~$y_0$ and of the length of the planning horizon~$T$. By defining the average expected charging rate and the average desired charging rate~as 
\begin{equation*}
    \dot y(x^b, x^r) = \frac{\E \left[ y(x^b, x^r, \tilde \delta, y_0, T) \right] - y_0}{T} = \eta^+ x^b - \eta_d x^r \varphi\left( - \frac{x^b}{x^r}\right)
    \quad \mathrm{and} \quad
    \dot y^\star = \frac{y^\star - y_0}{T},
\end{equation*}
respectively, the constraint~$\E [y(x^b, x^r, \tilde \delta, y_0, T)] = y^\star$ can be reformulated equivalently as~$\dot y(x^b, x^r) = \dot y^\star$. Since $T>0$, the expected charging rate~$\dot y$ inherits the concavity and monotonicity properties of the state-of-charge established in Proposition~\ref{prop:esoc}. In particular, if $x^r = 0$, then $\dot y(x^b, 0) = \eta^+ [x^b]^+ - \frac{1}{\eta^-} [x^b]^-$. Hence, $\dot y(x^b,0) = \dot y^\star$ is valid if and only if $x^b = \frac{1}{\eta^+} [\dot y^\star]^+ - \eta^- [\dot y^\star]^-$, which is fully determined by the desired charging rate~$\dot y^\star$ and by the charging and discharging efficiencies $\eta^+$ and~$\eta^-$. As $x^r$ increases, the expected charging rate~$\dot y(x^b, x^r)$ may decrease due to increased charging and discharging losses. The battery operator, however, can compensate this decrease by increasing~$x^b$. Since~$\dot y$ is strictly increasing, continuous, and surjective onto~$\mathbb{R}$ for any fixed~$x^r$, there is a unique~$x^b$ that satisfies the equation $\dot y(x^b, x^r) = \dot y^\star$. This means that~$x^b$ can be expressed as an implicit function of~$x^r$. This implicit function depends on the desired charging rate~$\dot y^\star$, the charging and discharging efficiencies, and the mean absolute deviation of the frequency deviations. The latter is defined as $\Delta = \mathbb{E}[\vert \tilde \xi \vert] = 2\varphi(0)$, where the second equality can be proved via integration by parts.

\begin{Rmk}[Mean absolute deviation]
\label{rmk:mad}
    We have~$\Delta \leq \frac{\gamma}{T}$ because all relevant frequency deviation trajectories reside in~$\D$ and have thus a mean absolute deviation no greater than~$\frac{\gamma}{T}$. Formally,
    \begin{equation*}
    \pushQED{\qed}
        \mathbb{P}\left[\delta \in \D\right] = 1 
        ~~ \implies ~~
        \Delta \leq \frac{\gamma}{T}.
        \qedhere
        \popQED
    \end{equation*}
\end{Rmk}

\begin{Prop}[Implicit function]\label{prop:if}
	The constraint $\dot y(x^b, x^r) = \dot y^\star$ defines a unique implicit function $g:\R_+ \to \R$ such that $\dot y(g(x^r), x^r) = \dot y^\star$ for all $x^r \in \mathbb{R}_+$. The function $g$ is convex, continuous, and nondecreasing with derivative 
	\begin{equation}\label{eq:g'}
		g'(x^r) = \eta_d \frac{\varphi(-\frac{x^b}{x^r}) + \frac{x^b}{x^r} F(-\frac{x^b}{x^r})}{\eta^+ + \eta_d F(-\frac{x^b}{x^r})}
	\end{equation}
	(if it exists), where $x^b = g(x^r)$. We have~$g'(x^r) = 0$ for all~$x^r \in (0, \vert g(0) \vert)$. In addition, the asymptotic slope~$m = \lim_{x^r \to \infty} g'(x^r) \in [0,1)$ is the unique solution of the equation 
	\begin{equation}\label{eq:ass}
		m = (1 - \eta^+ \eta^-) \varphi(m).
	\end{equation}
\end{Prop}

Proposition~\ref{prop:if} implies that for any fixed~$x^r$, the battery operator must buy the amount $x^b = g(x^r)$ of power in order to meet the expected state-of-charge target. One can interpret~$g(0)$ as the amount of power needed to meet the target in the absence of frequency regulation. Accordingly, $g(x^r) - g(0)$ reflects the amount of power needed to compensate the charging and discharging losses due to frequency regulation. These losses vanish for~$x^r = 0$ and increase in~$x^r$ at a rate that is smaller than or equal to~$m$. The asymptotic slope~$m$ is of particular interest because it is an upper bound on the percentage of regulation power that the battery operator needs to purchase in order to cover the losses from providing frequency regulation.

\begin{lem}[Asymptotic slope]\label{lem:m}
    The asymptotic slope~$m$ is convex and nonincreasing in the roundtrip efficiency~$\eta^+\eta^-$ and nondecreasing in the mean absolute deviation~$\Delta$ of the frequency deviations.
\end{lem}

Proposition~\ref{prop:if} further reveals that~$g'_+(0) = 0$ whenever~$y_0 \neq y^\star$. Otherwise, we have~$g'_+(0) = m$. 

\begin{lem}[Linearity]\label{lem:lin}
    If $y_0 = y^\star$,
    then the function~$g$ is linear with slope~$m$.
\end{lem}

If $y_0 = y^\star$, the losses increase exactly at rate~$m$. Maybe surprisingly, the losses are thus smaller when the initial state-of-charge differs from the target. 

\begin{Rmk}[Computability]\label{rmk:computability}
    For any fixed~$x^r$, the implicit function $g$ can be evaluated as follows. If $x^r = 0$, then $\dot y(g(0),0) = \dot y^\star$ readily implies that
    $g(0) = \frac{1}{\eta^+} [\dot y^\star]^+ - \eta^- [\dot y^\star]^-$. If $x^r > 0$, we first compute~$m$ as the unique solution to equation~\eqref{eq:ass} by bisection on~$[0,1]$. Next, we compute $g(x^r)$ as the unique root of the function $\dot y(g(x^r), x^r) - \dot y^\star$ on the interval~$[g(0), g(0) + mx^r]$, again by bisection. Once $x^b = g(x^r)$ is known, we obtain $g'(x^r)$ from equation~\eqref{eq:g'}. \hfill $\Box$
\end{Rmk}

Using Propositions~\ref{prop:cr} and~\ref{prop:if}, we can now reformulate problem~\eqref{pb:P} as the one-dimensional deterministic optimization problem
\begin{equation}
	\tag{P}
	\label{pb:P''}
	\begin{array}{>{\displaystyle}c*2{>{\displaystyle}l}}
		\min_{x^r \in \R_+} & \multicolumn{2}{>{\displaystyle}l}{
			T\left(c^b g(x^r) - c^r x^r\right)
		} \\
		\subj & x^r + g(x^r)   & \leq \bar{y}^+ \\
		& x^r - g(x^r)         & \leq \bar{y}^- \\
		& x^r + \max\left\{ \frac{T}{\gamma} g(x^r), g(x^r) \right\}  & \leq \frac{\bar{y} - y_0}{\eta^+ \gamma} \\
		& x^r - \min\left\{ \frac{T}{\gamma} g(x^r), g(x^r) \right\}  & \leq \frac{\eta^- y_0}{\gamma}. 
	\end{array}
\end{equation}

\begin{Th}[Constraint and dimensionality reduction]\label{th:dr}
    The problems~\eqref{pb:P} and~\eqref{pb:P''} are equivalent.
\end{Th}

Note that the objective function of problem~\eqref{pb:P''} is convex as the implicit function~$g$ is convex. The feasible set of~\eqref{pb:P''} can be represented concisely as $\set{X} = \left\{
    x^r \in \mathbb{R}_+ : \ell(x^r) \leq g(x^r) \leq u(x^r)
    \right\}$,~where
\begin{equation*}
    \ell(x^r) = \max\left\{
    x^r - \min\left\{ \bar y^-, \frac{\eta^- y_0}{\gamma}
    \right\}, \:
    \frac{\gamma}{T} x^r - \frac{\eta^- y_0}{T}
    \right\}
\end{equation*}
and
\begin{equation*}
    u(x^r) = \min\left\{
    \min\left\{
    \bar y^+, \frac{\bar y - y_0}{\eta^+ \gamma}\right\} - x^r, \:
    \frac{\bar y - y_0}{\eta^+ T} - \frac{\gamma}{T} x^r
    \right\}.
\end{equation*}

Due to the lower bounds on the convex function $g(x^r)$, the set~$\set{X}$ is nonconvex in general. However, if \mbox{$g(x^r) - \ell(x^r)$} is monotonic, then there exists at most one intersection between~$\ell$ and $g$, which means that the constraint $\ell(x^r) \leq g(x^r)$ defines nevertheless a convex feasible set.
This is the case under the following assumption.

\begin{Ass}[Roundtrip efficiency]\label{ass:ncnd_mad}
We have $\eta^+\eta^->\frac{1}{3}$.
\end{Ass}

Assumption~\ref{ass:ncnd_mad} is non-restrictive. Indeed, all relevant electricity storage technologies, as identified by the~\cite{WEC20}, have a roundtrip efficiency higher than~$\frac{1}{3}$.

\begin{lem}[Convex feasible set]\label{lem:cfs}
If Assumption~\ref{ass:ncnd_mad} holds, then the set~$\set{X}$ is convex.
\end{lem}

Lemma~\ref{lem:cfs} asserts that the feasible set~$\set{X}$ is a line segment under realistic parameter settings. As the objective function of problem~\eqref{pb:P''} is convex, an optimal solution~$x^r_\ast$ coincides either with a boundary point of the line segment~$\set{X}$ or with a stationary point of the objective function in the interior of the line segment. All three candidate solutions can be computed conveniently via bisection. If Assumption~\ref{ass:ncnd_mad} fails to hold, then~$\set{X}$ consists of two disjoint line segments, and it becomes necessary to check five different candidate solutions, all of which can be computed by bisection. A detailed discussion of this generalized setting is omitted because it has little practical relevance.
\section{Candidate Solutions}\label{sec:cs}
In the following, we will first examine the boundary points of the line segment~$\set{X}$ and then the stationary points of the objective function.

One can show that~\mbox{$u(x^r) - g(x^r)$} is strictly decreasing and that \mbox{$g(x^r) - \ell(x^r)$} is strictly decreasing under Assumption~\ref{ass:ncnd_mad}. The line segment~$\set{X}$ is thus non-empty if and only if $u(0) - g(0) \geq 0$ and $g(0) - \ell(0) \geq 0$. In this case, there exists a unique~$x^r_u$ such that $u(x^r_u) = g(x^r_u)$ and a unique~$x^r_\ell$ such that~$g(x^r_\ell) = \ell(x^r_\ell)$. The line segment~$\set{X}$ is thus given by \mbox{$\{ x^r \in \mathbb{R}_+ : x^r \leq x^r_\ell, \, x^r \leq x^r_u \} = [0, \bar x^r]$}, where $\bar x^r = \min\{ x^r_\ell, x^r_u\}$. At least one of the points $x^r_\ell$ and $x^r_u$ will be in~$[0, \bar{\bar x}^r]$, where $\bar{\bar x}^r$ is the unique intersection between the strictly increasing lower bound~$\ell$ and the strictly decreasing upper bound~$u$, and can be computed by bisection. The point~$\bar{\bar x}^r$ itself admits the closed form expression
 \begin{align*}
		\bar{\bar x}^r = \min \Bigg\{ & \frac{\bar y^+ + \bar y^-}{2}, \, 
		\frac{\bar y^+ + \frac{\eta^-}{\gamma} y_0}{2}, \,
		\frac{T \bar y^+ + \eta^- y_0}{\gamma + T}, \,
		\frac{\bar y^- + \frac{\bar y - y_0}{\eta^+ \gamma }}{2}, \,
		\frac{T\bar y^- + \frac{\bar y - y_0}{\eta^+}}{\gamma + T}, \\
		&
		\frac{\bar y + \big( \eta^+ \eta^- \frac{T}{\gamma} - 1 \big) y_0}{\eta^+(\gamma + T)}, \,
		\frac{ \frac{T}{\gamma}\bar y - \big(\frac{T}{\gamma} - \eta^+ \eta^-\big) y_0 }{\eta^+(\gamma + T)}
		\Bigg\}. 
\end{align*}

The stationary points of the objective function of problem~\eqref{pb:P''} are such that the expected marginal cost $Tc^b g'(x^r)$ of providing frequency regulation equals the expected marginal revenue $Tc^r$ from providing frequency regulation. The set of all stationary points in~$\set{X}$ is thus~$\set{X}_\star = \{ x^r \in [0,\bar x^r]: \frac{c^r}{c^b} \in \partial g(x^r) \}$, where $\partial g(x^r)$ denotes the subdifferential of $g$ at $x^r$. Note that~$\frac{c^r}{c^b} > 0$ because~$c^b > 0$ and~$c^r > 0$. If~$\bar x^r = 0$, then~$0$ is the only feasible solution to problem~\eqref{pb:P''}. If~$g'_-(\bar x^r) < \frac{c^r}{c^b}$,
then~$\set{X}_\star$ is empty, and~$\bar x^r$ is the optimal solution to problem~\eqref{pb:P''} because the marginal revenue of providing frequency regulation is strictly higher than the marginal cost for all feasible~$x^r$. Similarly, if~$g'_+(0) > \frac{c^r}{c^b}$, then~$\set{X}_\star$ is again empty, and~$0$ is the optimal solution because the marginal cost of providing frequency regulation is strictly higher than the marginal revenue for all~$x^r \in (0, \infty]$. Otherwise,~$\set{X}_\star$ is non-empty and may contain several stationary points, all of which are optimal solutions to problem~\eqref{pb:P''}. In this case, we assume that the battery operator selects the smallest stationary point to avoid unnecessary battery usage. Theorem~\ref{th:x_ast} formalizes these results.

\begin{Th}\label{th:x_ast}
For $\bar x^r \geq 0$, the smallest optimal solution to problem~\eqref{pb:P''} is
\begin{equation*}
    x^r_\ast = \begin{cases}
     0 & \text{if } \bar x^r = 0 \text{ or } \bar x^r > 0 \text{ and } g'_+(0) \geq \frac{c^r}{c^b}, \\
    \bar x^r & \text{if } \bar x^r > 0 \text{ and } g'_-(\bar x^r) < \frac{c^r}{c^b}, \\
    \min \set{X}_\star & \text{otherwise}.
    \end{cases}
\end{equation*}
\end{Th}

\begin{Rmk}
If~$\dot y^\star = 0$, then $g(x^r) = m x^r$. Therefore, $x^r_\ast = \bar x^r$ if $m < \frac{c^r}{c^b}$ and $x^r_\ast = 0$ otherwise. \hfill $\Box$
\end{Rmk}

\begin{Rmk}\label{rmk:inf_X}
If $\set{X}_\star$ is non-empty, $\min \set{X}_\star $ can be found by bisection on~$\set{X}$ as $g'$ is nondecreasing. \hfill $\Box$
\end{Rmk}

To compute the optimal solution $x^r_\ast$, we first compute~$\bar x^r$. If~$\bar x^r > 0$, then we evaluate~$g'(0)$ and~$\min \partial g(\bar x^r)$. By Proposition~\ref{prop:if} and Lemma~\ref{lem:lin}, we find~$g'(0) = 0$ if~$\dot y^\star \neq 0$; $ = m$ otherwise. Finally, if~$\set{X}_\star$ is non-empty, we compute its minimum. Since $\bar x^r$, $g'_-(\bar x^r)$, and $\min \set{X}_\star$ can all be computed by bisection, the optimal solution~$x^r_\ast$ can also be computed highly efficiently by bisection.

\begin{figure}
    \centering
    \includegraphics[width = \textwidth]{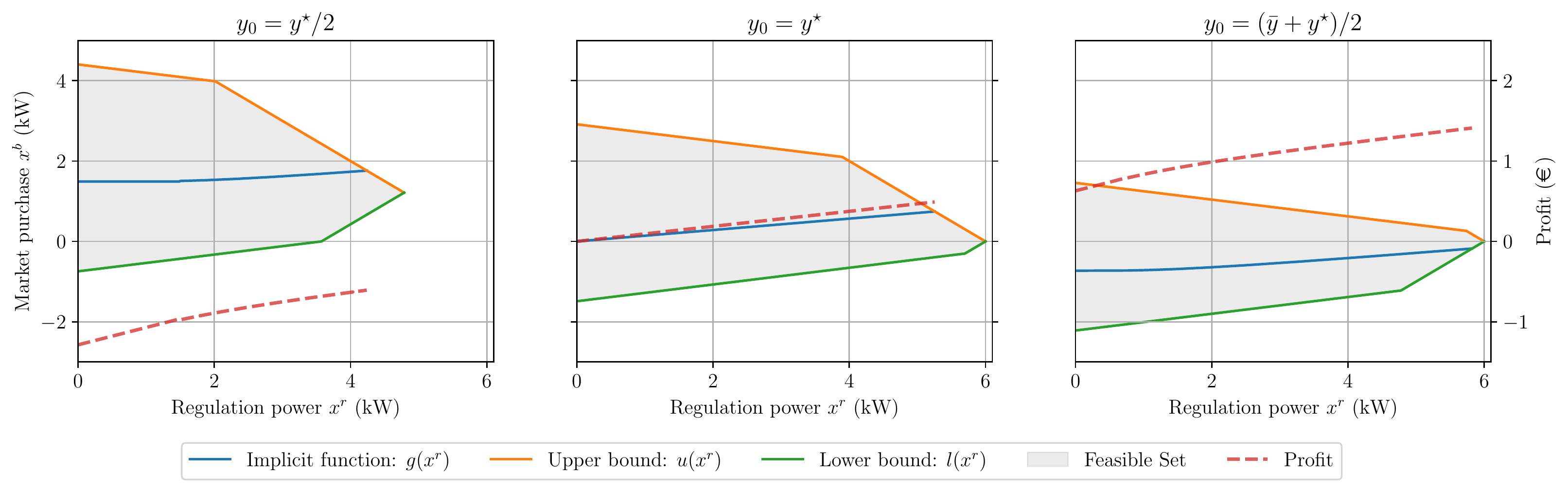}
    \caption{Feasible set and profit for $\eta^+ = \eta^- = 0.707$, $\Delta = 5 \cdot 0.0816$, $\gamma = 5$h, $T = 24$h, $\bar y = 100$kWh, $\bar y^+ = 6$kW, $\bar y^- = 6$kW, $c^r = 0.9$cts/kWh, $c^{b} = 4 c^r$.}
    \label{fig:X}
\end{figure}

Figure~\ref{fig:X} illustrates the dependence of the feasible set on the initial state-of-charge. In the particular example depicted here, the net profit is increasing in $x^r$ for all $x^r \in [0, \bar x^r]$. It is thus optimal to set $x^r$ to $\bar x^r$, which is given by the intersection of the graph of the implicit function $g$ with the graphs of the lower bound $l$ or the upper bound $u$.

The expected marginal cost of providing frequency regulation depends on the desired charging rate~$\dot y^\star$, which is only known to the battery operator but unknown to the grid operator. Nevertheless, the grid operator knows that the expected marginal cost amounts to at most~$c^b m$ and can therefore infer that it is profitable for the battery operator to offer all available regulation power if~$\frac{c^r}{c^b} > m$. From Lemma~\ref{lem:m}, we know that~$m$ is convex and nonincreasing in the roundtrip efficiency and nondecreasing in the mean absolute deviation of the frequency deviations. In the next section, we will derive explicit lower and upper bounds on~$m$ that are tight for certain degenerate frequency deviation distributions. For these distributions, we will derive analytical solutions to problem~\eqref{pb:P''}.

\begin{Rmk}
    Section~\ref{sec:si_elastic_prices} in the Supplementary Material extends Theorem~\ref{th:x_ast} to elastic prices. \hfill $\Box$
\end{Rmk}

\section{Analytical Solution}\label{sec:as}
We now construct two discrete distributions~$\ubar{ \mathbb{P}}_\xi$ and~$\bar{\mathbb{P}}_\xi$ with the same mean absolute deviation as~$\mathbb{P}_\xi$, which is given by~$\Delta = 2\varphi(0)$. Specifically, we define~$\ubar{\mathbb{P}}_\xi$ as a two-point distribution with mass~$\frac{1}{2}$ at~$-\Delta$ and~$\Delta$, and $\bar{\mathbb{P}}_\xi$ as a three-point distribution with mass~$\frac{\Delta}{2}$ at~$-1$ and~$1$, and mass $1-\Delta$ at~$0$. The super-cumulative distribution functions~$\ubar \varphi$ and $\bar \varphi$ of~$\ubar{\mathbb{P}}_\xi$ and~$\bar{\mathbb{P}}_\xi$ with~$\ubar \varphi(-1) = 0$ and $\bar \varphi(-1) = 0$ are ~$\ubar \varphi(\xi) = \max\{ 0, \frac{1}{2}( \xi + \Delta), \xi \}$ and $\bar \varphi(\xi) = \max\{ 0, \frac{\Delta}{2}(\xi + 1), (1-\frac{\Delta}{2})\xi + \frac{\Delta}{2}, \xi \}$, respectively.

\begin{lem}[Second-order stochastic dominance]\label{lem:phi}
    We have~$\ubar \varphi(\xi) \leq \varphi(\xi) \leq \bar \varphi(\xi)$ for all~$\xi \in \mathbb{R}$.
\end{lem}

We now define the asymptotic sensitivities~$\underline m$ and~$\overline m$ as the solutions to the nonlinear algebraic equations \mbox{$\underline m = (1 - \eta^+ \eta^-)\ubar \varphi(\underline m)$} and $\overline m = (1 - \eta^+\eta^-) \bar \varphi(\overline m)$, respectively, which exist and are unique by Proposition~\ref{prop:if}. These equations admit the closed-form solutions
\begin{equation*}
    \underline m =  \frac{1-\eta^+\eta^-}{1+\eta^+\eta^-}\Delta
    \quad
    \text{and}
    \quad
    \overline m = 1 - \frac{1}{1 + (\frac{1}{\eta^+\eta^-}-1)\frac{\Delta}{2}}.
\end{equation*}

The discrete distributions~$\ubar{\mathbb{P}}_\xi$ and~$\bar{\mathbb{P}}_\xi$ provide explicit bounds on the implicit function~$g(x^r)$. These bounds, which we denote by~$\ubar g(x^r)$ and~$\bar g(x^r)$, are obtained by solving the differential equation~\eqref{eq:g'} for~$\ubar{\mathbb P}_\xi$ and~$\bar{\mathbb P}_\xi$, respectively. Since $\ubar \varphi$ and~$\bar \varphi$ are piecewise linear, the differential equations can be solved separately for each linear piece. Combining the results for the different pieces yields
\begin{align*}
	\ubar g(x^r) = & \max\left\{
	g(0),~
	\underline m x^r + g(0) - \frac{1 - \eta^+ \eta^-}{1+\eta^+\eta^-} \vert g(0) \vert
	\right\}\\ \text{and} \quad
	\bar g(x^r) = & \max\left\{
	g(0),~
	\frac{(1-\eta^+\eta^-)\varphi(0)x^r -[g(0)]^-}{\eta^+\eta^-+(1-\eta^+\eta^-)(1-\varphi(0))},~
	\overline m x^r + \frac{\eta^+\eta^- [g(0)]^+ - [g(0)]^-}{\eta^+\eta^-+(1-\eta^+\eta^-)\varphi(0)}
	\right\}.
\end{align*}
\begin{lem}\label{lem:g_bds}
    We have~$\ubar g(x^r) \leq g(x^r) \leq \bar g(x^r)$ for all~$x^r \in \mathbb{R}_+$ and \mbox{$\underline m \leq m \leq \overline m$}.
\end{lem}

In order to calculate the optimal solution $x^r_\ast$ under~$\ubar{\mathbb{P}}_\xi$ and~$\bar{\mathbb{P}}_\xi$, we need to compute the smallest stationary points, if they exist, as well as the boundary points of the feasible sets~$\ubar{\set{X}}$ and~$\bar{\set{X}}$. If they exist, the smallest stationary points occur either at kinks of the piecewise linear objective functions $c^b \ubar g(x^r) - c^r x^r$ and $c^b \bar g(x^r) - c^r x^r$ or at~$x^r = 0$ because~$\ubar g$ and~$\bar g$ are piecewise linear.
The left boundary points of both feasible sets are~0, as under~$\mathbb{P}_\xi$. The right boundary points are~$\min\{\ubar x^r_\ell, \bar x^r_u\}$ and $\min\{ \bar x^r_\ell, \ubar x^r_u \}$, respectively, where $\ubar x^r_\ell$, $\bar x^r_u$, $\bar x^r_\ell$, and~$\ubar x^r_u$ are the unique points with $\ubar g(\ubar x^r_\ell) = \ell(\ubar x^r_\ell)$, $\ubar g(\bar x^r_u) = u(\bar x^r_u)$, 
$\bar g(\bar x^r_\ell) = \ell(\bar x^r_\ell)$, 
and $\bar g(\ubar x^r_u) = u(\ubar x^r_u)$. The points $\ubar x^r_\ell$, $\bar x^r_u$, $\bar x^r_\ell$, and~$\ubar x^r_u$ can be computed in closed form by evaluating the intersections of the affine functions~$\ell$ and~$u$ with the affine functions that generate~$\ubar g$ and~$\bar g$. The right boundary points $\min\{\ubar x^r_\ell, \bar x^r_u\}$ and $\min\{ \bar x^r_\ell, \ubar x^r_u \}$ are then given by the minimum of six and eight rational functions of the problem parameters~$y_0$, $\dot y^\star$, $\bar y$, $\bar y^+$, $\bar y^-$, $\eta^+$, $\eta^-$, $\Delta$, $\gamma$, and~$T$, respectively. We omit the closed-form expressions because they are too cumbersome to be insightful in general. 

In order to derive insightful solutions, we assume that $y_0 = y^\star$ in the subsequent analysis. In this special case, the losses due to frequency regulation are higher than for any other value of~$y_0$ as explained in the discussion after Proposition~\ref{prop:if}. If it is profitable to provide frequency regulation in this special case, it will thus also be profitable to provide frequency regulation in any other case. As~$\dot y^\star = \frac{y^\star - y_0}{T} = 0$, Lemma~\ref{lem:lin} implies that $\ubar g(x^r) = \underline{m} x^r$, $g(x^r) = mx^r$, and~$\bar g(x^r) = \overline{m} x^r$. Under any frequency deviation distribution~$\mathbb{P}_\xi$, the optimal solution to problem~\eqref{pb:P''} then satisfies~$x^r_\ast = \bar x^r$ if $m < \frac{c^r}{c^b}$ and~$x^r_\ast = 0$ otherwise. Since~$\ell$ and~$u$ are piecewise linear functions with two pieces each, the right boundary point~$\bar x^r$ of the feasible set can be expressed as the minimum of only four rational functions of the problem parameters, and will admit an intuitive interpretation.

\begin{Prop}\label{prop:xr_lin}
If~$y_0 = y^\star$, then
\begin{equation}\label{eq:xr}
    \bar x^r = \min\left\{
    \frac{\bar y^-}{1-m},~
    \frac{\bar y^+}{1+m},~
    \frac{\eta^- y_0}{\gamma(1-m)},~
    \frac{\bar y - y_0}{\eta^+(\gamma + m T)}
    \right\}.
\end{equation}
\end{Prop}
The first two terms in formula~\eqref{eq:xr} for~$x^r$ depend on the charging capacity~$\bar y^+$ and the discharging capacity~$\bar y^-$, while the last two terms depend on the battery capacity~$\bar y$ and the initial state-of-charge~$y_0$. We say that the battery is \emph{power-constrained} if~$\bar x^r$ is equal to one of the first two terms. Otherwise, we say that the battery is~\emph{energy-constrained}. We now state the analytical solution.

\begin{Th}[Analytical solution]\label{th:as}
If~$y_0 = y^\star$, then an optimal solution to problem~\eqref{pb:P''} is
\begin{equation*}
    x^r_\ast = \begin{cases}
    0 & \text{if } m \geq \frac{c^r}{c^b}, \\
    \bar x^r & \text{otherwise},
    \end{cases}
\end{equation*}
under any frequency deviation distribution~$\mathbb{P}_\xi$. If~$\mathbb{P}_\xi = \ubar{\mathbb{P}}_\xi$, then $m = \underline{m}$.  If~$\mathbb{P}_\xi = \bar{\mathbb{P}}_\xi$, then $m = \overline{m}$.
\end{Th}

As a direct consequence of Theorem~\ref{th:as},
if $\overline{m} < \frac{c^r}{c^b}$, it is optimal to set~$x^r = \bar x^r$ for any frequency deviation distribution with mean absolute deviation~$\Delta$, regardless of the shape of the distribution, because~$\overline{m} \geq m$. In the following, we describe the maximum amount~$\bar x^r$ of regulation power that can be offered by power- and energy-constrained batteries for the case~$y_0 = y^\star$.

If the battery is power-constrained, then~$\bar x^r$ depends on the charging and discharging efficiencies only through the marginal increase $m$ in the expected power loss which, in turn, depends on these efficiencies only through the roundtrip efficiency~$\eta^+\eta^-$. Due to charging and discharging losses, the battery operator expects to lose energy while delivering frequency regulation and compensates the expected loss by purchasing the power~$mx^r$ from an electricity market, which decreases the effective charging capacity and increases the effective discharging capacity of the battery. Accounting for the effective charging and discharging capacities, the battery operator may dimension the discharging capacity~$\bar y^-$ as a fraction~$\frac{1-m}{1+m}$ of the charging capacity~$\bar y^+$ without restricting~$\bar x^r$.

If the battery is energy-constrained, then~$\bar x^r$ depends not only on the roundtrip efficiency~$\eta^+\eta^-$, through~$m$, but also on the individual charging and discharging efficiencies. Charging losses increase the amount of energy that the battery can consume from the grid and therefore increase the effective storage capacity. Conversely, discharging losses decrease the amount of energy that the battery can deliver to the grid and therefore decrease the effective storage capacity. For given charging and discharging efficiencies, the initial state-of-charge~$y_0$ determines how much energy the battery can consume from and deliver to the grid. As the battery operator must be able to both consume and deliver regulation power, $\bar x^r$~is maximized if the battery can absorb as much energy from the grid as it can deliver to the grid. This occurs at an initial state-of-charge of~$y^\star_0$ and results in the maximum amount~$\bar x^r_\star$ of regulation power that can be offered by energy-constrained batteries, where
\begin{equation*}
    \bar x^r_\star = \min\left\{
    \frac{\bar y^-}{1-m},~
    \frac{\bar y^+}{1+m},~
    \frac{\eta^-}{\frac{\gamma}{T}(1 + \eta^+\eta^- - m) + \eta^+\eta^-m} \frac{\bar y}{T}
    \right\}
    \text{ and }
    y^\star_0 = \frac{(1 - m) \bar y}{1 + \eta^+\eta^- + (\eta^+\eta^- \frac{T}{\gamma} - 1)m}.
\end{equation*}
Assuming that~$\bar y^- = \frac{1-m}{1+m} \bar y^+$, the battery is thus energy constrained if $\frac{\bar y^+}{1+m} > \bar x^r_\star $, which occurs if the battery's charge rate $C = \frac{\bar y^+}{y}$ (C-rate) is no smaller than
\begin{equation*}
    \underline{C} = \frac{1}{T} \cdot \frac{(1+m)\eta^-}{\frac{\gamma}{T}(1+\eta^+\eta^--m)+\eta^+\eta^-m}.
\end{equation*}
The C-rate expresses the percentage of the battery's storage capacity that can be consumed from the grid within one hour. 

The initial state-of-charge~$y^\star_0$, which maximizes the amount of regulation power that energy-constrained batteries can provide, depends on the charging and discharging efficiencies only through the roundtrip efficiency~$\eta^+\eta^-$. It increases from~$\frac{\bar y}{2}$ to~$\bar y$ as the roundtrip efficiency decreases from~$1$ to~$0$. The maximum amount~$\bar x^r_\star$ of regulation power that can be provided by energy-constrained batteries is equal to~$\frac{\bar y}{2\gamma}$ in the absence of charging and discharging losses. The storage capacity~$\bar y$ is divided by~$2\gamma$ because the battery operator must be able to both consume and deliver all of the regulation power she promised for a total time of at least~$\gamma$. 
\section{Case Studies}\label{sec:appl}
In the following, we first parameterize our model based on the French frequency regulation market. Next, we analyze the marginal cost and profit of providing frequency regulation as well as the maximum amount of regulation power that storage operators can provide. Last, we discuss the profits that storage operators can earn over the planning horizon~$\T$ and, for the specific case of lithium-ion batteries, under what conditions on the uncertainty budget~$\gamma$ and the length of the planning horizon~$T$ it may be profitable to invest in electricity storage for frequency~regulation. 

\subsection{Model Parametrization}\label{sec:param}
We focus on storage operators who provide frequency regulation to the French grid operator and compute their profits based on historical frequency deviation data, on availability and delivery prices, and on wholesale and retail market prices. Frequency measurements, availability prices, and delivery prices are published by the French grid operator Réseau de Transport d'Electricité (RTE).\footnote{\url{https://clients.rte-france.com}} Wholesale market prices depend on how long before delivery electricity is traded. Since frequency regulation is traded up to one day before delivery, we use the prices of the day-ahead market, which are equal to the delivery prices published by RTE \cite[p.~68]{RTE17}. Retail market prices vary from one electric utility company to another. We use the base tariff of Electricité de France, the largest French electric utility company. The French government regulates this particular tariff and publishes the corresponding prices.\footnote{Journal Officiel de la République Fran\c{c}aise: \url{https://legifrance.gouv.fr}} We will compare the operating costs of storage technologies with different roundtrip efficiencies, namely hydrogen, redox flow batteries, vehicle-to-grid, pumped hydro, and stationary lithium-ion batteries. In assigning roundtrip efficiencies to storage technologies, we follow~\cite{EA17} for vehicle-to-grid and the World Energy Council \cite[p.~9]{WEC20} for all other storage technologies. Table~\ref{tab:effs} lists the roundtrip efficiencies of the storage technologies. In order to judge whether it is profitable to invest in lithium-ion batteries for frequency regulation, we follow~\cite{SC19} and assume that investment costs in the year~$2023$ range from US\$$85$ to US\$$165$ per~kWh of storage capacity with a lifetime of $10$~years and from US\$$710$ to US\$$860$ per kW of charging and discharging capacity with a lifetime of $30$~years. We assume an exchange rate of~$1\text{\EUR{}} = \text{US}\$1.15$ and annualize the investment costs with a yearly discount rate of~$2\%$, which equals the long-term inflation target of the~\cite{ECB21}. We use 2019 electricity prices, \ie, before the drop in prices during the Covid pandemic and the rise in prices since the war in Ukraine and the maintenance problems in the French nuclear power plant fleet. Based on these prices, we set the ratio $\frac{c^r}{c^b}$ of the expected average price of regulation power to the expected market price of electricity to~$ 0.251$ for wholesale market prices and to~$0.059$ for retail market prices. In terms of frequency regulation, we consider the common European market for frequency containment reserves, which has a daily planning horizon and thus set~$T  = 24$~hours. We approximate the cumulative distribution function corresponding to~$\mathbb{P}_\xi$ by the symmetric logistic function~$F(\xi) = \frac{1}{1 + \exp(-\theta \xi)}$ with $\theta = \frac{2 \ln(2)}{\Delta}$ and mean absolute deviation $\Delta = 0.0816$. Finally, we set $\dot y^\star = 0$ in all experiments. SM \ref{sm:params} provides more details about the model parametrization. We provide all code and data at \url{www.github.com/lauinger/cost-of-frequency-regulation-through-electricity-storage} and on Zenodo \citep{data}.

\begin{table}[t!]
    \centering
    \caption{Roundtrip efficiencies of storage technologies.}
    \begin{tabular}{c|ccccc}
        Storage Technology & Hydrogen & Redox Flow & Vehicle-to-Grid & Pumped Hydro & Li-Ion \\
        \midrule
        Roundtrip Efficiency (\%) & 35--55 & 60--85 & 70--85 & 75--85 & 85--95
    \end{tabular}
    \label{tab:effs}
\end{table}

\subsection{Marginal Cost and Maximum Regulation Bid}
\label{sec:m_xr}
\begin{figure}[t!]
	\centering
	\includegraphics[width=1\linewidth]{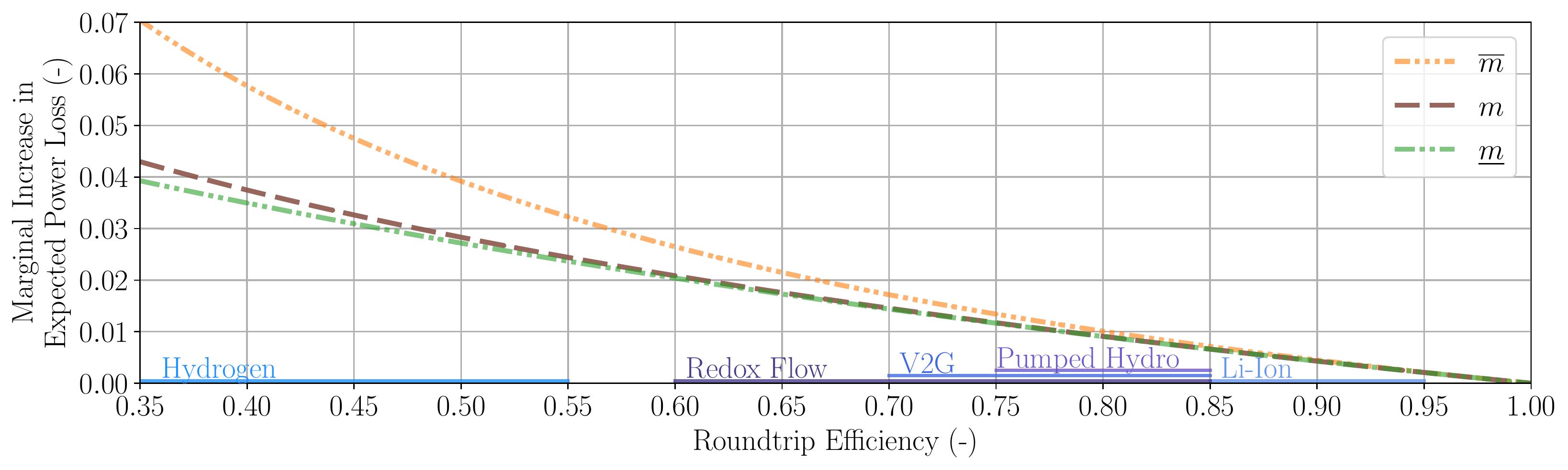}
	\caption{Marginal increase in expected power loss for different roundtrip efficiencies.}
	\label{fig:m}
\end{figure}

The marginal cost~$mc^bT$ of providing frequency regulation is the marginal increase in the expected power loss~$m$ multiplied by the market price of electricity~$c^b$ and the length of the planning horizon~$T$. It is profitable to use a storage device at its full potential if~$mc^bT$ is lower than the marginal revenue $Tc^r$, \ie, if and only if the ratio~$\frac{c^r}{c^b}$ of the expected average price of regulation power to the expected average market price of electricity exceeds the marginal increase~$m$ in the expected power loss. Figure~\ref{fig:m} displays~$m$ as a function of the roundtrip efficiency for the estimated logistic distribution of frequency deviations, together with its lower bound~$\underline{m}$ and its upper bound~$\overline{m}$, both parametrized to have the same mean absolute deviation as the logistic distribution. The bounds are tight when the roundtrip efficiency equals one and loosen as the roundtrip efficiency decreases. The upper bound loosens faster than the lower bound. For roundtrip efficiencies higher than~$0.60$, the lower bound, $\underline{m} = \frac{1-\eta^+\eta^-}{1+\eta^+\eta^-}\Delta$, underestimates~$m$ by less than~$4.59\cdot 10^{-4}$. At a roundtrip efficiency of~$0.35$, typical for inefficient hydrogen storage, $m$ equals~$4.30\%$. For inefficient redox flow batteries, with a roundtrip efficiency of~$0.60$, $m$ decreases to~$2.09\%$. For inefficient lithium-ion batteries, with a roundtrip efficiency of~$0.85$, $m$ decreases further to~$0.66\%$. Unsurprisingly, $m$ vanishes for perfectly efficient storage devices. For storage operators buying electricity at retail prices, the ratio~$\frac{c^r}{c^b}$ was greater than~$0.026$ on all days in~2019. It would have therefore been profitable for them to use any of the storage technologies we consider at their full potential for frequency regulation, except for hydrogen storage. For storage operators buying electricity at wholesale prices, the ratio~$\frac{c^r}{c^b}$ was greater than~$0.07$ on all days in~2019, which means that it would have been profitable for them to use a hydrogen storage device at its full potential, too. 

\begin{figure}[t!]
	\centering
	\includegraphics[width=1\linewidth]{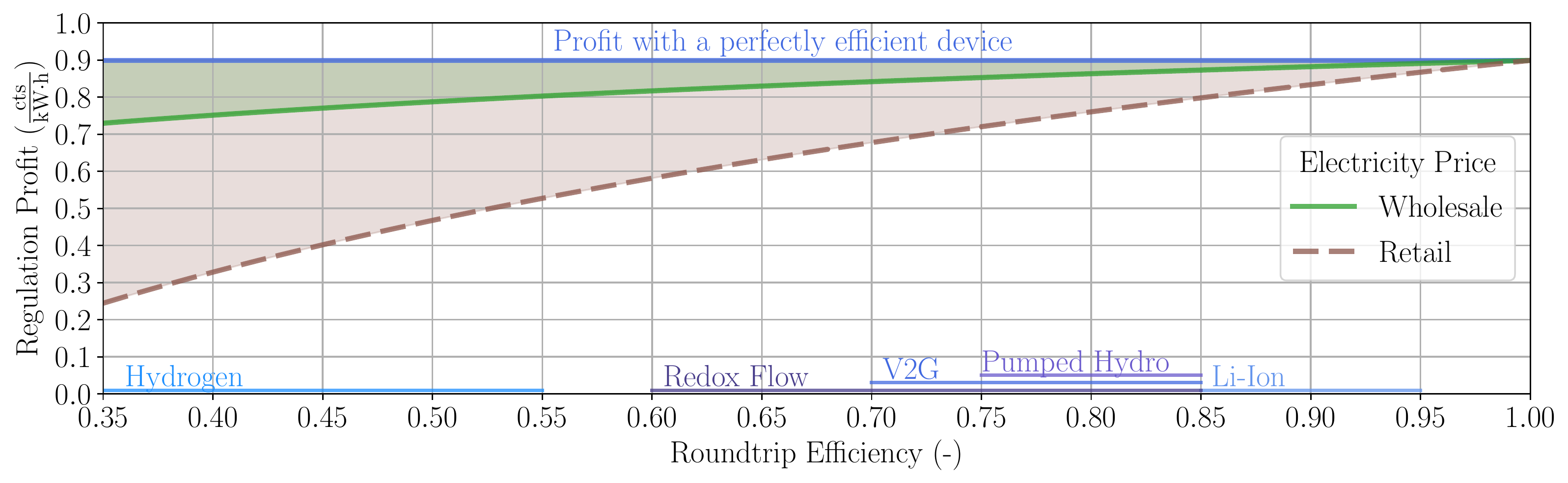}
	\caption{Profit per unit of regulation power for wholesale and retail electricity prices.}
	\label{fig:profit_kW}
\end{figure}

Although charging and discharging losses may not cause storage operators to withhold regulation power from the market, they still reduce the profit~$c^r - mc^b$ per unit of regulation power. Figure~\ref{fig:profit_kW} shows the average profit per unit of regulation power as a function of roundtrip efficiency for storage operators buying electricity at wholesale and retail prices in the year~2019. Losses reduce the profit by~$19\%$ from $0.90\frac{\text{cts}}{\text{kW} \cdot \text{h}}$ to $0.73\frac{\text{cts}}{\text{kW} \cdot \text{h}}$ for the most inefficient storage devices if electricity is bought at wholesale prices. If electricity is bought at retail prices, losses reduce the profit by~$72\%$ to $0.24\frac{\text{cts}}{\text{kW} \cdot \text{h}}$. The reduction is four times higher at retail than at wholesale prices. A hydrogen tank with~$ \eta^-\eta^+ = 0.4$ buying electricity at wholesale prices and an electric vehicle with~$\eta^-\eta^+ = 0.8$ buying electricity at retail prices achieve the same profit per unit of regulation power. 

Losses reduce not only the profit per unit of regulation power, but may also impact the amount of regulation power~$\bar x^r$ a storage device can provide. If the storage device is energy-constrained, then charging losses increase the amount of energy that the storage device can consume from the grid, while discharging losses decrease the amount of energy that the storage device can deliver to the grid. In principle, charging losses may outweigh discharging losses and increase the normalized regulation bid compared to storage devices with no losses. In practice, however, discharging losses usually outweigh charging losses. As examples, we consider lithium-ion batteries, vehicle-to-grid, and hydrogen storage, all of them operating at an \emph{activation ratio}~$\frac{\gamma}{T}$ of $0.2$, which measures the fraction of time during which a storage operator must be able to provide all the regulation power she promised.
Hydrogen storage is unlikely to be energy-constrained because hydrogen can be stored at low cost in steel tanks, or at even lower cost in salt caverns~\citep{MV19}. Nevertheless, we include hydrogen in our comparison as an example of storage devices with low roundtrip efficiencies. For lithium-ion batteries with charging and discharging efficiencies of~$0.92$ each, losses reduce the normalized regulation bid from~$1$ to~$0.98$. For vehicle-to-grid with a charging efficiency of~$0.88$ and a discharging efficiency of~$0.79$, loosely based on~\cite{EA17}, losses reduce the normalized regulation bid to~$0.91$. For hydrogen storage with a charging efficiency of~$0.80$ and a discharging efficiency of~$0.58$, based on~\cite{MV19}, losses reduce the normalized regulation bid to~$0.77$. For more details on the impact of charging and discharging losses on~$\bar x^r$, see SM \ref{sm:xr}.

\subsection{Profits}
We will now analyze the profit that an energy-constrained storage device may earn per unit of storage capacity. First, we describe the \emph{operating} profit made over the planning horizon of length~$T$. For lithium-ion batteries, we then calculate the \emph{effective yearly} profit as the difference between the operating profits made over one year and the annualized investments costs of the battery.

\begin{figure}[!t]
\centering
\begin{tikzpicture}
\tikzstyle{every node}=[font=\small]
    \begin{axis}[
        legend style = {
        cells = {anchor=west},
        anchor = north east,
        },
        ytick = {nc, nd, cts-kW-ws, cts-kW-rt, max-reg, cts-kWh-ws, cts-kWh-rt},
        xtick distance = 0.25,
        y tick label style  = {align=right, rotate = 0},
        x tick label style={
        /pgf/number format/.cd,
            fixed,
            fixed zerofill,
            precision=2,
        /tikz/.cd
        },
        grid = both,
        xmin = 0, xmax = 2.275,
        width = 12cm,
        height = 0.325\linewidth,
        symbolic y coords = {cts-kWh-rt, cts-kWh-ws, max-reg, cts-kW-rt, cts-kW-ws, nd, nc},
        yticklabels = {%
        Charging Efficiency $(-)$,
        Discharging Efficiency $(-)$, 
        Wholesale Regulation Profit $(\frac{\text{cts}}{\text{kW$\cdot$ h}})$,
        Retail Regulation Profit $(\frac{\text{cts}}{\text{kW$\cdot$ h}})$,
        Normalized Regulation Power $(-)$,
        Wholesale Operating Profit $(\frac{\text{cts}}{\text{kWh}})$,
        Retail Operating Profit $(\frac{\text{cts}}{\text{kWh}})$
        },
        ]
    \addlegendimage{empty legend}
    \addlegendentry{\hspace{-0.8cm} Storage Technology}
    \addplot[xbar, draw = black, fill = black!50, opacity = 0.5] table [y=label, x = Ideal] {comparison.txt};
    \addlegendentry{Ideal Storage}
    \addplot[xbar, draw = black, fill = CornflowerBlue!40, opacity = 1] table [y=label, x = Lithium-Ion] {comparison.txt};
    \addlegendentry{Li-Ion Battery}
    \addplot[xbar, draw = black, fill = RoyalBlue!60, opacity = 1] table [y=label, x = V2G] {comparison.txt};
    \addlegendentry{Vehicle-to-Grid}
    \addplot[xbar, draw = black, fill = dodgerblue!100, opacity = 1] table [y=label, x = H2] {comparison.txt};
    \addlegendentry{Hydrogen}
    \end{axis}
\end{tikzpicture}
\caption{The impact of charging and discharging losses on regulation power and profits.}
\label{fig:losses}
\end{figure}
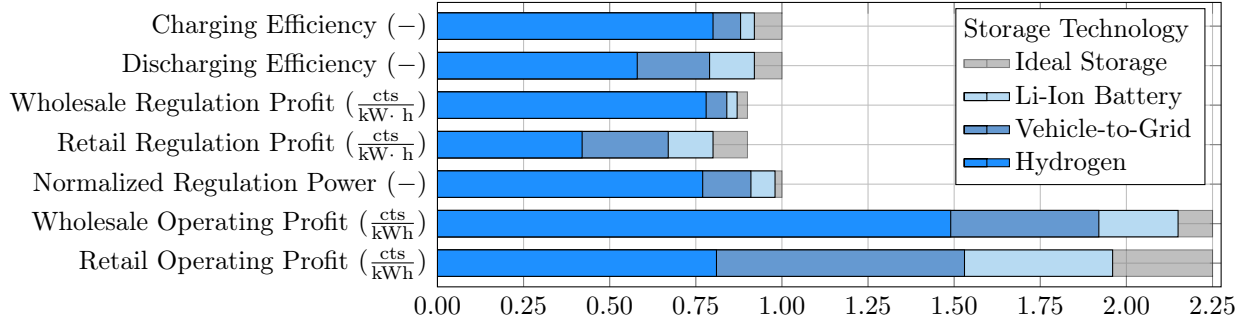

The operating profit is the product of the profit per unit of regulation power and the amount of regulation power the storage device can deliver. Formally, the operating profit is thus equal to
\begin{equation*}
    (c^r - m c^b) \cdot \frac{\eta^- \bar y}{\frac{\gamma}{T}(1 + \eta^+\eta^--m) + \eta^+\eta^-m}.
\end{equation*}
It may be somewhat surprising that the operating profit depends on the length of the planning horizon~$T$ only through the activation ratio~$\gamma / T$. Intuitively, one may expect the operating profit to increase linearly with the length of the planning horizon because a longer planning horizon should allow the storage operator to sell the same amount of regulation power for a longer period of time. This reasoning is correct if the storage device is power-constrained. For energy-constrained storage devices, however, the storage operator can only deliver a fixed amount of regulation energy. The length of the planning horizon does therefore not influence the operating profit. Figure~\ref{fig:losses} summarizes the impact of losses on the profits per unit of regulation power, on the maximum normalized regulation power, and on operating profits at wholesale and retail electricity prices for the lithium-ion battery, vehicle-to-grid, and hydrogen storage examples that we considered earlier. SM \ref{sm:op} explains in detail how operating profits depends on the charging and discharging efficiencies.

\begin{figure}[t!]
	\centering
	\includegraphics[width=\linewidth]{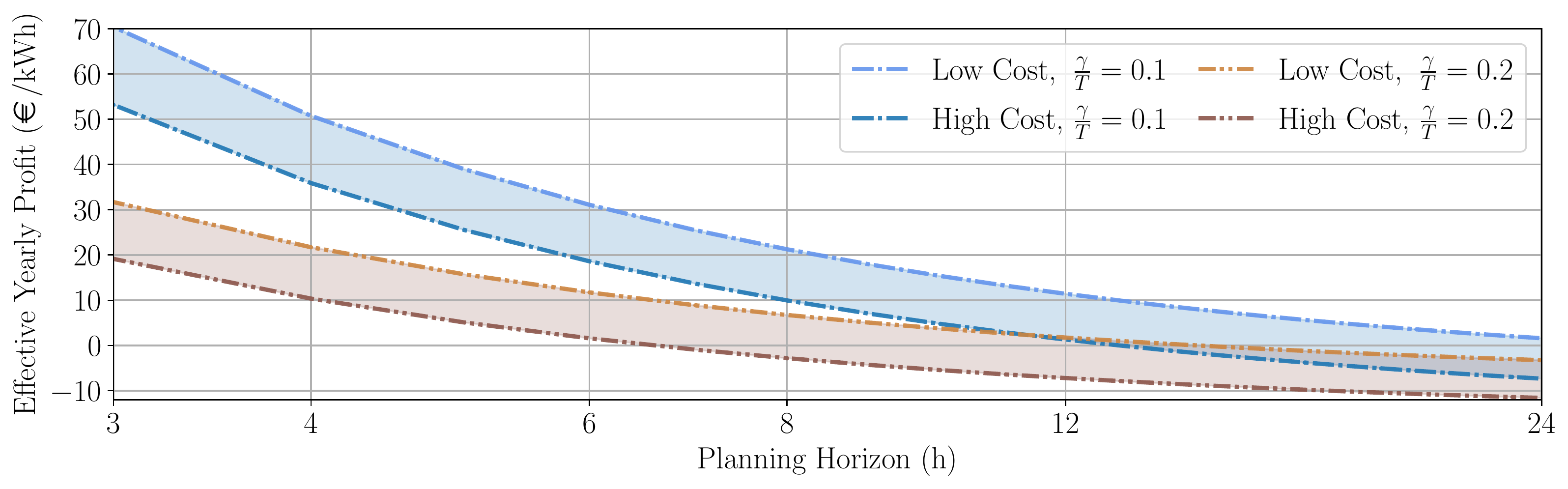}
	\caption{Effective yearly profit per kWh of storage capacity after investment costs.}
	\label{fig:pi_y}
\end{figure}

The effective yearly profit determines whether it is worthwhile for a storage operator to invest into storage devices for frequency regulation. Based on the cost and lifetime data of lithium-ion batteries in Section~\ref{sec:param}, we estimate the annualized costs of lithium-ion batteries in the year 2023 to range from $8.2\frac{\text{\EUR{}}}{\text{kWh}}$ to $16.0\frac{\text{\EUR{}}}{\text{kWh}}$ for energy storage capacity and from $27.6\frac{\text{\EUR{}}}{\text{kW}}$  to $33.4\frac{\text{\EUR{}}}{\text{kW}}$ for charging and discharging capacity. 
Figure~\ref{fig:pi_y} shows the effective yearly profit for lithium-ion batteries with charging and discharging efficiencies of~$0.92$, buying electricity at wholesale prices, as a function of the length of the planning horizon, for activation ratios of~$0.1$ and~$0.2$, for low annualized investment costs of $8.2\frac{\text{\EUR{}}}{\text{kWh}}$ and $27.6\frac{\text{\EUR{}}}{\text{kW}}$, and for high annualized investment costs of $16.0\frac{\text{\EUR{}}}{\text{kWh}}$ and  $33.4\frac{\text{\EUR{}}}{\text{kW}}$. At the current $24$~hour planning horizon, lithium-ion batteries are profitable only at an activation ratio of~$0.1$ and low investment costs. Given that we only considered the cost of the battery itself but no additional costs related to installation, maintenance, administration, or land lease, investing in lithium-ion batteries for frequency regulation does not seem to be profitable in the near future. In the medium term, lithium-ion batteries may become sufficiently low-cost~\citep{MZ21} to be used for frequency regulation. Besides falling battery prices, grid operators might opt for an activation ratio of~$0.1$ rather than~$0.2$ to make the use of energy storage for frequency regulation more profitable. This would roughly double the operating profits from frequency regulation, but it would also shrink the uncertainty set~$\D$ and therefore make grid operators more vulnerable to extreme frequency deviations that could cause black-outs.

Alternatively, grid operators could reduce the length of the planning horizon~$T$. If the planning horizon were to be reduced from 24~hours to 4~hours, for example, the operating profits accrued over a one-year period would increase by a factor~$6$. Similarly, the minimum C-rate required for the battery to be energy-constrained and thus the costs for charging and discharging capacity would also increase by a factor~$6$. The increase in operating profits is well worth the increased charger costs. Under the reduced planning horizon, lithium-ion batteries could achieve an effective yearly profit of~$10$\EUR per kWh of storage capacity, even at high investment costs and high activation ratios. They would achieve this higher profit by selling more regulation power, which may decrease prices for frequency regulation and reduce marginal profits. Section~\ref{sec:si_elastic_prices} in the supplementary material shows that the price effect is negligible for an individual 100kWh lithium-ion battery. It may, however, become more significant if many battery operators sell more regulation power. A shorter planning horizon for frequency regulation does not necessarily make grid operators more vulnerable to extreme frequency deviations and is already common practice in intraday wholesale electricity markets.\footnote{\url{https://www.entsoe.eu/network_codes/cacm/implementation/sidc/}} Adopting intraday markets for frequency regulation as well would make energy-constrained storage devices more competitive with power-constrained flexibility providers, such as thermal power plants and pumped hydro storage. The increased competition may decrease the total cost of frequency regulation, which is ultimately borne by the public.

\section{Conclusions}

The \cite{ENTSOE22} estimates that they will need to store electricity in batteries with a total power of up to~$240$GW by the year~2050. Lithium-ion batteries are considered a promising source of frequency regulation, thanks to their fast dynamics. The investment costs of lithium-ion batteries have declined sharply in recent years, but we find that they are not yet low enough for lithium-ion batteries to be profitable in the frequency regulation market. Since Europe plans to increasingly rely on batteries in the future, their use should become profitable. 

We identify two policy options that make electricity storage in general and battery storage in particular more profitable in the frequency regulation market. First, regulators can decrease the marginal costs of frequency regulation by making it easier for small and medium-sized storage devices to access wholesale electricity markets. This is one of the aims of Order~845 by the US \cite{FERC18}. Second, regulators can decrease the length of the planning horizon, which is currently one day in the common European frequency regulation market. We show that the amount of regulation power that storage devices can provide may be constrained by their storage capacity and their initial state-of-charge. In this case, the profits from frequency regulation over the lifetime of the storage devices are inversely proportional to the length of the planning horizon. The planning horizon could be shortened by adopting intraday markets for frequency regulation, which already exist for the wholesale of electricity. 

\paragraph{Acknowledgements} D.L. thanks Emilia Suomalainen, Ja\^{a}far Berrada, Fran\c{c}ois Colet, Willett Kempton, and Yannick Perez for helpful discussions, and the Institut Vedecom for funding.

\linespread{1}
\small
\bibliographystyle{abbrvnat}
\bibliography{bibliography}
\linespread{1.5}
\normalsize
\appendix
\section{Proofs}\label{sec:ch3_apx_proofs}
This appendix contains the proofs of all essential theorems, propositions, and lemmas in the main text. All other proofs are relegated to the supplementary material. The proof of Proposition~\ref{prop:cr} relies on the following lemma, which involves the set~$\D_\downarrow$ of all nonincreasing left-continuous functions in~$\D$.

\begin{lem}\label{lem:D_down}
    For any $x^r \geq 0$, we have
    \begin{equation}\label{eq:D_down}
        \max_{\delta \in \D, \, t \in \T} y(x^b, x^r, \delta, y_0, t)
        =
        \max_{\delta \in \D^+_\downarrow, \, t \in \T} y(x^b, x^r, \delta, y_0, t).
    \end{equation}
\end{lem}

\begin{proof}[Proof of Proposition~\ref{prop:cr}]
    We first show that the upper bound on the charging power and the upper bound on the state-of-charge are valid for all frequency deviation trajectories~$\delta \in \D$ and all time instants~$t \in \T$ if and only if they are valid for the particular frequency deviation trajectory~$\delta^{(+)}$, defined through $\delta^{(+)}(t) = 1$ if $t \leq \gamma$ and $\delta^{(+)} = 0$ otherwise, and both time instants~$t \in \{\gamma,T\}$.
    
    The upper bound on the charging power is valid for all~$\delta \in \D$ and all~$t \in \T$ if and only if it is valid for the maximum charging power that can be achieved by any~$\delta \in \D$ and any~$t \in \T$. We have
        \begin{equation*}
        \max_{\delta \in \D,\,t \in \T} y^+(x^b, x^r, \delta(t)) 
        = \max_{\delta \in \D^+, \, t \in \T} y^+(x^b, x^r, \delta(t))
        = \max_{\delta \in \D^+, \, t \in \T} x^b + \delta(t) x^r 
        = x^b + x^r,
    \end{equation*}
    where the first equality holds because~$y^+$ is nondecreasing in~$\delta(t)$ and because~$\D$ is symmetric. In fact, for any~$\delta \in \D$, we have $\vert \delta\vert \in \D$, and the maximum charging power for~$\vert \delta \vert$ is at least as high as the one for~$\delta$. The second equality holds because~$y^+$ is linear in~$\delta(t)$ whenever~$\delta(t) \geq 0$. The last equality holds because $\delta(t) \leq 1$ for all~$\delta \in \set{D}^+$ and~$t \in \set{T}$, and because the upper bound is attained at~$\delta = \delta^{(+)}$ and~$t = \gamma$, for example. Thus, assertion~$(i)$ follows. 
    
    For the upper bound on the state-of-charge, we first use Lemma~\eqref{lem:D_down} to obtain
    \begin{equation*}
        \max_{\delta \in \D, \, t \in \T} y(x^b, x^r, \delta, y_0, t)
        =
        \max_{\delta \in \D^+_\downarrow, \, t \in \T} y(x^b, x^r, \delta, y_0, t).
    \end{equation*}
    If $x^b + x^r < 0$, then the battery is discharging for all $t \in \T$. The upper bound on the state-of-charge is thus valid if $y_0 \leq \bar y$, which we stated as a condition in Proposition~\ref{prop:cr}.

    Otherwise, if $x^b + x^r \geq 0$, one can show that
    \begin{equation}\label{eq:cr_2d}
        \max_{\delta \in \D^+_\downarrow, \, t \in \T} y(x^b, x^r, \delta, y_0, t)
        = \left\{ \arraycolsep=4pt\def\arraystretch{1.2}
	    \begin{array}{>{\displaystyle}c*1{>{\displaystyle}l}}
		\max_{t^c, \delta^c} 
		& y_0 + t^c \eta^+ (x^b + \delta^c x^r) \\
		\subj & x^b + \delta^c x^r \geq 0,~~
		t^c \delta^c \leq \gamma,~~0 \leq \delta^c \leq 1,~~0 \leq t^c \leq T. 
	\end{array}
    \right.
    \end{equation}
    To this end, we first show that any feasible solution to the maximization problem on the left-hand side of~\eqref{eq:cr_2d} gives rise to a feasible solution to the maximization problem on the right-hand side with the same or a larger objective value. For any $\delta \in \D^+_\downarrow$ and $t \in \T$, we construct the last time at which the battery is still charged as $t^c = \max_{t' \in [0,t]} \{t': x^b + \delta(t') x^r \geq 0\}$, which exists because~$\delta$ is left-continuous nonincreasing and hence upper semi-continuous. We also construct the average frequency deviation signal during the charging process as $\delta^c = \frac{1}{t^c} \int_0^{t^c} \delta(t') \, \text{d}t'$, which satisfies $x^b + \delta^c x^r \geq 0$ because $x^b + \delta(t') x^r \geq 0$ for all~$t \leq t^c$ because $\delta$ is nonincreasing. As $\delta \in \D^+_\downarrow$, we have $t^c \delta^c = \int_0^{t^c} \delta(t') \, \text{d}t' \leq \int_\T \delta(t') \, \text{d}t' \leq \gamma$. In addition, we have $0 \leq t^c \leq T$ as $t \in \T$ and $0 \leq \delta^c \leq 1$ as $\delta \in  \D^+_\downarrow$. Since the state-of-charge $y$ is nondecreasing in~$\delta$, which is nonincreasing in time, we have
    \begin{equation}\label{eq:cr_2d_obj}
        y(x^b, x^r, \delta, y_0, t)
        \leq
        y(x^b, x^r, \delta, y_0, t^c)
        =
        \int_0^{t^c }y_0 + \eta^+ (x^b + \delta(t')x^r) \, \text{d}t' 
        =
        y_0 + t^c \eta^+(x^b + \delta^c x^r),
    \end{equation}
    where the second equality holds because $\delta$ is integrated against a constant function and can thus be set to its average value~$\delta^c$ over the integration horizon.

    Given a feasible solution~$(t^c, \delta^c)$ to the maximization problem on the right-hand side of~\eqref{eq:cr_2d}, we can also construct a feasible solution to the left-hand side with the same objective value by setting~$t = t^c$ and $\delta(t') = \delta^c$ if $t' \leq t^c$ and $=0$ otherwise. It is clear that $t \in \T$ since $0 \leq t^c \leq T$. In addition, $\delta$ is nonnegative and left-continuous nonincreasing as $0 \leq \delta^c \leq 1$. Finally, $\int_\T \delta(t') \, \text{d}t' = t^c \delta^c \leq \gamma$ ensures that $\delta \in \D^+_\downarrow$. The construction satisfies again $t^c = \max_{t' \in [0,t]} \{t': x^b + \delta(t') x^r \geq 0\}$ and $\delta^c = \frac{1}{t^c} \int_0^{t^c} \delta(t') \, \text{d}t'$. The equality of the objective values thus follows from~\eqref{eq:cr_2d_obj}.
  
    In summary, we have shown that the two problems in~\eqref{eq:cr_2d} have indeed the same maximum.
    
    The objective function of the right-hand side problem in~\eqref{eq:cr_2d} is nondecreasing in $\delta^c$ and $t^c$ since $\eta^+ \geq 0$, $x^r\geq 0$, and $x^b + \delta^c x^r \geq 0$. It is thus optimal to make $\delta^c t^c$ as large as possible, \ie, setting it to~$\gamma$, because the box constraints on~$\delta^c$ and $t^c$ only imply a weaker upper bound of $T \geq \gamma$ on $t^c \delta^c$. By substituting $\delta^c = \gamma / t^c$ we arrive at the equality
    \begin{equation}\label{eq:cr_1d}
    \left. \arraycolsep=1pt\def\arraystretch{1.2}
    \begin{array}{>{\displaystyle}c*1{>{\displaystyle}l}}
		\max_{t^c, \delta^c} 
		& y_0 + t^c \eta^+ (x^b + \delta^c x^r) \\
		\subj & x^b + \delta^c x^r \geq 0, \,
		t^c \delta^c = \gamma, \,
		0 \leq \delta^c \leq 1, \,
        0 \leq t^c \leq T 
	\end{array}
    \right\}
    = 
    \left\{ \arraycolsep=1pt\def\arraystretch{1.2}
	   \begin{array}{>{\displaystyle}c*1{>{\displaystyle}l}}
		\max_{t^c} 
		& y_0 + \eta^+ (t^c x^b + \gamma x^r) \\
		\subj & t^c x^b + \gamma x^r \geq 0, \,
        \gamma \leq t^c \leq T.
	\end{array}
    \right.
    \end{equation}
    We will now analyze the one-dimensional linear program on the right-hand side of~\eqref{eq:cr_1d}. If $x^b \geq 0$, then the inequality $t^c x^b + \gamma x^r \geq 0$ is always valid and it is optimal to make $t^c$ as large as possible, \ie, setting it to~$T$. Conversely, if $x^b < 0$, then it is optimal to make~$t^c$ as small as possible, \ie, setting it to~$\gamma$. In this case $t^c$ satisfies again the inequality $t^c x^b + \gamma x^r \geq 0$ because~$x^b + x^r \geq 0$. Combining all of the above arguments, we thus obtain
    \begin{equation*}
        \max_{\delta \in \D, \, t \in \T} y(x^b, x^r, \delta, y_0, t)
        =
        y_0 + \eta^+ \big( \max\{ \gamma x^b, T x^b \} + \gamma x^r \big),
    \end{equation*}
    leading to assertion~(\emph{iii}). One easily verifies that $\delta^{(+)}$ attains the maximum on the left-hand side.
  
  Using similar arguments, one can show that the upper bound on the discharging power and the lower bound on the state-of-charge hold for all frequency deviation signals~$\delta \in \D$ and all time instants~$t \in \T$ if and only if they hold for the particular frequency deviation signal~$\delta^{(-)} = - \delta^{(+)}$ and all time instants~$t\in\{\gamma, T\}$. We omit the details for the sake of brevity.
  \end{proof}

The proof of Proposition~\ref{prop:esoc} relies on the following symmetry property of~$\varphi$.
\begin{lem}[Symmetry of~$\varphi$]\label{lem:sym}
	For all $z \in \R$, we have $\varphi(z) = \varphi(-z) + z$.
\end{lem}

\begin{proof}[Proof of Proposition~\ref{prop:esoc}]
    We first prove equation~\eqref{eq:esoc} for $x^r > 0$. We have
    \begin{align}
        \E \left[
        y(x^b, x^r, \tilde \delta, y_0, T)
        \right]
        & =
        y_0 + T \E \left[
        \frac{1}{T}\int_\T \eta^+ \left[ x^b + \tilde \delta(t) x^r\right]^+ 
        - \frac{1}{\eta^-} \left[
        x^b + \tilde \delta(t) x^r
        \right]^-
        \, \mathrm{d}t
        \right] \notag \\ \label{eq:esoc_pr1}
        & = y_0 + T x^r \int_{-1}^1 \eta^+ \left[ \frac{x^b}{x^r} + \xi \right]^+ - \frac{1}{\eta^-} \left[ \frac{x^b}{x^r} + \xi \right]^- \, \mathbb{P}_\xi(\mathrm{d}\xi),
    \end{align}
    where the second equality follows from the definition of~$\mathbb{P}_\xi$. Setting $z = x^b/x^r$ to simplify notation, 
	\begin{align*}
		\int_{-1}^{1} [ z + \xi ]^+ \, \mathbb{P}_\xi(\mathrm{d} \xi) &
		= \int_{-z}^{1} ( z + \xi ) \, \mathbb{P}_\xi(\mathrm{d} \xi)
		= z F(\xi) \rvert_{-z}^1 + \xi F(\xi) \rvert_{-z}^1 - \int_{-z}^{1} F(\xi) \, \mathrm{d} \xi \\
	    & = (z + 1) F(1) + \varphi(-z) - \varphi(1)
		= z + \varphi(-z) .
	\end{align*}
	The second equality follows from integration by parts, whereas the fourth equality holds because $F(1) = \varphi(1) = 1$. In fact, as $\varphi(-1) = 0$ by construction, Lemma~\ref{lem:sym} implies that 
	$\varphi(1) = \varphi(-1) +1 = 1$.
	Following a similar reasoning and keeping in mind that $F(-1) = 0$, we obtain
	\begin{align*}
		\int_{-1}^{1} [ z + \xi ]^- \, \mathbb{P}_\xi(\mathrm{d} \xi) &
		= -\int_{-1}^{-z} ( z + \xi ) \, \mathbb{P}_\xi(\mathrm{d} \xi)
		= (z - 1) F(-1) + \varphi(-z) - \varphi(-1)  = \varphi(-z).
	\end{align*}
    Substituting these expressions into~\eqref{eq:esoc_pr1} yields equation~\eqref{eq:esoc}.
 
    For~$x^r = 0$, the expectation $\mathbb{E}[y(x^b,0, \tilde \delta, y_0, T)]$ is trivially equal to $y_0 + T (\eta^+ [x^b]^+ - \frac{1}{\eta^-} [x^b]^-)$, which corresponds to the formula given in Proposition~\ref{prop:esoc} thanks to the definition of the perspective function. In fact,     
    \begin{equation*}
        \lim_{x^r \to 0^+} x^r \varphi\left( - \frac{x^b}{x^r} \right)
        = x^b \left( 
        \lim_{x^r \to 0^+} \frac{\partial}{\partial x^b} \, x^r \varphi\left(-\frac{x^b}{x^r}\right)
        \right)
        = x^b \left(
        \lim_{x^r \to 0^+} - F\left( -\frac{x^b}{x^r} \right)
        \right)
        = \left[ x^b \right]^-.
    \end{equation*}

    We now establish several useful properties of the expected terminal state-of-charge. Note first that~$\varphi$ is convex, Lipschitz continuous, and almost everywhere differentiable because it is a super-cumulative distribution function. The expected terminal state-of-charge is jointly concave in~$x^b$ and~$x^r$ because~$-x^r\varphi(-\frac{x^b}{x^r})$ is the negative perspective of the convex function~$\varphi(-x^b)$ and therefore concave \citep[p.~89]{SB04}. As a perspective of a Lipschitz continuous convex function, the expected terminal state-of-charge is also globally continuous.     To see that the expected terminal state-of-charge is strictly increasing in~$x^b$ for every $x^r > 0$, note that
    \begin{equation*}
        \frac{\partial}{\partial x^b} \E \left[
        y(x^b, x^r, \tilde \delta, y_0, T)
        \right]
        = T\left( \eta^+ + \eta_d F\left( - \frac{x^b}{x^r} \right) \right) > 0
        \quad \forall (x^b, x^r) \in \mathbb{R} \times \mathbb{R}_{++}
    \end{equation*}
   because $\eta^+ > 0$, $\eta_d \geq 0$, and $F$ is nonnegative. For~$x^r = 0$, $\mathbb{E}[y(x^b,0, \tilde \delta, y_0, T)] = y_0 + T (\eta^+ [x^b]^+ - \frac{1}{\eta^-} [x^b]^-)$, which is strictly increasing in~$x^b$ since $\eta^+ > 0$ and $\eta^- > 0$.

   Similarly, to see that 
   the expected terminal state-of-charge is nondecreasing in~$x^r$, we note that
    \begin{equation*}
        \frac{\partial}{\partial x^r} \E \left[
        y(x^b, x^r, \tilde \delta, y_0, T)
        \right]
        =
        - \eta_d T\left( \varphi\left(-\frac{x^b}{x^r}\right) + \frac{x^b}{x^r}F\left(-\frac{x^b}{x^r}\right)\right)
        \leq 0
        \quad \forall (x^b, x^r) \in \mathbb{R} \times \mathbb{R}_{++}.
    \end{equation*}
    To prove the inequality, we set $z = -x^b / x^r$ and show that the function $-\eta_d (\varphi(z) - z F(z))$ is nonnegative. As $\varphi(z) = 0$ for all $z \leq -1$, Lemma~\ref{lem:sym} implies that $\varphi(z) = z$ for all $z \geq 1$. Thus, we have $\varphi(z) - zF(z) = 0$ for all $\vert z \vert \geq 1$. If $z \in [-1, 0]$, then $\varphi(z) - zF(z) \geq 0$ because $\varphi$ and $F$ are both nonnegative. Finally, if~$z \in [0,1]$ then we first note that
    \begin{equation}\label{eq:esoc_pr2}
        \varphi(1) = \varphi(z) + \int_{z}^1 F(z') \, \text{d}z' 
        \leq \varphi(z) + \int_z^1 F(1) \, \text{d}z'
        = \varphi(z) + F(1)(1-z).
    \end{equation}
    Hence, we have
    \begin{equation*}
        0 = \varphi(1) - F(1)
        \leq \varphi(z) + F(1)(1-z) - F(1)
        = \varphi(z) - z F(1)
        \leq \varphi(z) - z F(z)
    \end{equation*}
    for every~$z \in [0,1]$, where the second inequality follows from the monotonicity of~$F$. 
    
    Finally, the expected terminal state-of-charge is unbounded above in~$x^b$ because
    \begin{equation*}
        \lim_{x^b \to \infty} y_0 + T\left(
        \eta^+ x^b - \eta_d x^r \varphi\left( -\frac{x^b}{x^r} \right)
        \right)
        = \lim_{x^b \to \infty} y_0 + T \eta^+ x^b
        = \infty,
    \end{equation*}
    where the first equality holds because~$\varphi(-\frac{x^b}{x^r}) = 0$ for all $x^b \geq x^r$.
\end{proof}

\begin{proof}[Proof of Proposition~\ref{prop:if}]
    The average expected charging rate~$\dot y$ is a positive affine transformation of the expected terminal state-of-charge. Proposition~\ref{prop:esoc} thus immediately implies that $\dot y$ is continuous and jointly concave in~$x^b$ and~$x^r$. In addition, $\dot y$ is strictly increasing and unbounded above in~$x^b$, and nonincreasing in~$x^r$. As~$\dot y$ is concave and strictly increasing in~$x^b$, it is also unbounded below in~$x^b$. Overall, $\dot y$ is continuous and unbounded below and above in~$x^b$, which means that the equation~$\dot y(x^b,x^r) = \dot y^\star$ has at least one solution~$x^b$ for any given~$x^r \in \mathbb{R}_+$. As~$\dot y$ is strictly increasing in~$x^b$, this solution is also unique. The constraint~$\dot y(x^b, x^r) = \dot y^\star$ defines therefore a unique implicit function $g:\mathbb{R}_+ \to \mathbb{R}$ such that $\dot y (g(x^r), x^r) = \dot y^\star$ for all~$x^r \in \mathbb{R}_+$. 
    
    As~$\dot y(x^b, x^r)$ is nonincreasing in~$x^r$ and strictly increasing in~$x^b$, the equality~$\dot y(g(x^r), x^r) = \dot y^\star$ remains valid if and only if the implicit function~$g$ is nondecreasing.
    
    As~$\dot y(x^b, x^r)$ is concave in~$x^b$ and~$x^r$, the superlevel set \mbox{$\set{C} = \{ (x^b, x^r) \in \mathbb{R} \times \mathbb{R}_+: \dot y(x^b, x^r) \geq \dot y^\star\}$} is convex. As~$\dot y$ is strictly increasing in~$x^b$, a point~$(x^b,x^r)$ satisfies $\dot y(x^b, x^r) \geq \dot y^\star$ if and only if $x^b \geq g(x^r)$. The set~$\set{C}$ thus coincides with the epigraph of~$g$. The convexity of~$\set{C}$ then implies that~$g$ is a convex function~\citep[p.~75]{SB04}.
     
    The proof of Proposition~\ref{prop:esoc} reveals that the partial derivatives of the expected terminal state of charge with respect to $x^b$ and $x^r$ exist on $\mathbb{R}$ and $\mathbb{R}_{++}$, which implies that the partial derivatives of~$\dot y$ also exist on $\mathbb{R}$ and $\mathbb{R}_{++}$. Since~$\dot y$ is partially differentiable and continuous, the univariate function~$g$ will be continuous and differentiable almost everywhere. By the implicit function theorem~\citep[p.~395]{MP85}, the derivative of~$g$ is given by
    \begin{equation*}
		g'(x^r)
		= -\frac{\frac{\partial \dot y(x^b,\, x^r)}{\partial x^r}}{\frac{\partial \dot y(x^b,\, x^r)}{\partial x^b}}
		= \eta_d \frac{\varphi(-\frac{x^b}{x^r}) + \frac{x^b}{x^r} F(-\frac{x^b}{x^r})}{\eta^+ + \eta_d F(-\frac{x^b}{x^r})}
	\end{equation*}
	if it exists. To show that $g'(x^r) = 0$ for every $x^r \in (0, \vert g(0) \vert)$, assume that $g(0) \neq 0$, and note that
	\begin{equation*}
	    \dot y(g(0), \vert g(0) \vert)
	    = \eta^+ g(0) - \eta_d \, \vert g(0) \vert \, \varphi\left( -\frac{g(0)}{\vert g(0)\vert} \right)
	    = \eta^+ g(0) - \eta_d \left[g(0)\right]^-
	    = \dot y(g(0), 0),
	\end{equation*}
	where the first and third equalities follow from Proposition~\ref{prop:esoc}, and the second equality holds because~$\varphi(z) = 0$ for~$z \leq -1$ and~$\varphi(z) = z$ for~$z \geq 1$. This implies that $\dot y (g(\vert g(0) \vert), \vert g(0) \vert) = \dot y(g(0), 0)$, and thus $g(\vert g(0) \vert) = g(0)$. As~$g$ is nondecreasing, it must be constant throughout the interval~$[0, \vert g(0) \vert]$, which means that~$g'(x^r) = 0$ for all~$x^r$ in the interior of that interval. 
	
	Note that the asymptotic slope of $g$ is given by~$m = \lim_{x^r \to \infty} g'(x^r) = \lim_{x^r \to \infty} \frac{g(x^r)}{x^r}$ because~$g$ is convex. The limit exists and is bounded below because $g$ is convex. In addition, we have
    \begin{equation*}
        0 = \lim_{x^r \to \infty} \frac{\dot y^\star}{x^r}
        = \lim_{x^r \to \infty} \frac{\dot y(g(x^r), x^r)}{x^r}
        = \lim_{x^r \to \infty} \eta^+ \frac{g(x^r)}{x^r} - \eta_d \varphi\left( - \frac{g(x^r)}{x^r}\right),
    \end{equation*}
    implying that $\lim_{x^r \to \infty} \frac{g(x^r)}{x^r} < +\infty$ as $\lim_{x^r \to \infty} \varphi(-\frac{g(x^r)}{x^r}) = 0$. Thus, $m = \lim_{x^r \to \infty} \frac{g(x^r)}{x^r}$ is finite, and~$m$ is a solution to the equation~$\eta^+ \mu - \eta_d \varphi(-\mu) = 0$. As~$\eta_d = \frac{1}{\eta^-} - \eta^+$ and \mbox{$\varphi(-m) = \varphi(m) - m$} by Lemma~\ref{lem:sym}, this equation is equivalent to~$m = (1 - \eta^+\eta^-)\varphi(m)$. It admits a unique solution within the interval~$[0,1)$, because the function $s(\mu) = \mu - (1 - \eta^+\eta^-)\varphi(\mu)$ is nondecreasing in~$\mu$, nonpositive for~$\mu = 0$, and strictly positive for~$\mu = 1$. The function $s$ is nondecreasing as its derivative \mbox{$1 - (1-\eta^+\eta^-)F(m)$} is nonnegative because~$F(m) \leq 1$ for all~$m \in \mathbb{R}$ and because~$\eta^+\eta^- \in (0,1]$. 
\end{proof}

\begin{proof}[Proof of Lemma~\ref{lem:m}]
    We know from Proposition~\ref{prop:if} that the asymptotic slope~$m$ is the unique solution to the equation $s(\mu) = 0$, where $s(\mu) = \mu - (1-\eta^+\eta^-)\varphi(\mu)$. We have $s(0) = -(1-\eta^+\eta^-)\varphi(0) \leq 0$ and~$s'(\mu) = 1 - (1-\eta^+\eta^-)F(\mu)$. As~$F(\mu) \in [0,1]$ for all~$\mu \in \mathbb{R}$ and as~$\eta^+\eta^-\in (0,1]$, $s'$ is nonnegative and $s$ is nondecreasing in~$\mu$. An increase in $\Delta = 2\varphi(0)$ can decrease (but not increase) the intercept $s(0)$, but it does not influence the slope~$s'$. An increase in~$\eta^+\eta^-$ can increase (but not decrease) the intercept $s(0)$ and the slope~$s'$. The zero-crossing of~$s$ is thus nondecreasing in~$\Delta $ and nonincreasing in~$\eta^+\eta^-$, and so is~$m$.

    To show that $m$ is convex and nonincreasing in the roundtrip efficiency $\eta^+\eta^-$, we first note that
    \begin{equation}
    \label{eq:m_leq}
     m = \min_\mu \left\{ \mu \, : \, \mu \geq (1-\eta^+\eta^-)\varphi(\mu)\right\}
    \end{equation}
    as it is never optimal to set $\mu > (1-\eta^+\eta^-)\varphi(\mu)$ since $\varphi$ is nondecreasing and continuous.
    
    Next, we characterize~$\varphi$ as a pointwise supremum of affine functions.
    As~$\varphi$ is a continuous convex function, it is closed. By the envelope representation theorem, any closed convex function is the pointwise supremum of all affine functions below it~\citep[p.~102]{RR70}. For~$\varphi$, specifically, we have~$\varphi(\mu) = 0$ for all~$\mu \leq -1$ and~$\varphi(\mu) = \mu$ for all~$\mu \geq 1$. It thus suffices to consider all affine functions~$a \mu + b$ with~$a \geq 0$ and~$b \geq 0$ such that $\varphi(\mu) \geq a \mu +b$. The highest possible slope of any such function is the highest possible slope of~$\varphi$, which is~$1$. Similarly, the highest possible intercept of any such function is the intercept of~$\varphi$, which is~$\varphi(0)$. Let~$\set{A} = \{ (a, b) \in [0, 1] \times [0, \varphi(0)] : \varphi(\mu) \geq a \mu + b~~\forall \mu \in \mathbb{R} \}$ be the set of all admissible coefficients for the affine functions. By the envelope representation theorem, we have~$\varphi(\mu) = \max_{(a,b) \in \set{A}} a \mu + b$ for all~$\mu \in \mathbb{R}$.    

    Substituting this expression of~$\varphi$ into equation~\eqref{eq:m_leq} yields
    \begin{align*}
        m = & \min_\mu \left\{ \mu \, : \, \mu \geq (1-\eta^+\eta^-) \max_{(a,b) \in \set{A}} a \mu + b \right\}
        = \min_\mu \left\{ \mu \, : \, \mu \geq (1-\eta^+\eta^-) a \mu + b  ~~ \forall (a,b) \in \set{A} \right\} \\
        = & \min_\mu \left\{ \mu \, : \, (1 - a(1-\eta^+\eta^-))\mu \geq (1 - \eta^+\eta^-)b  ~~ \forall (a,b) \in \set{A} \right\} \\
        = & \max_{(a,b) \in \set{A}} \varsigma(a,b,\eta^+\eta^-), 
        ~\text{where}~~
        \varsigma(a,b,\eta^+\eta^-) = \frac{1 - \eta^+ \eta^-}{1 - a (1 - \eta^+ \eta^-)}b.
    \end{align*}
    The first equality follows directly from substitution. The second equality holds because the inequality in the optimization problem is valid for all $(a,b) \in \set{A}$ if and only if it is valid for a pair $(a,b) \in \set{A}$ that maximizes the right-hand side of the inequality. Note that the embedded maximization in the inequality does indeed maximize the right-hand side since $1 - \eta^+\eta^- \geq 0$. The fourth equality holds because $1- a(1-\eta^+\eta^-) > 0$ as~$a \leq 1$ and~$\eta^+\eta^- \in (0,1]$. 

    We will now show that~$\varsigma$ is convex in~$\eta^+\eta^-$ for all~$(a,b) \in \set{A}$, which will later imply that~$m$ is also convex in~$\eta^+\eta^-$. Note that~$\varsigma$ is twice differentiable in~$\eta^+\eta^-$. In fact,
    \begin{equation*}
        \frac{\partial}{\partial (\eta^+\eta^-)} \,
        \varsigma(a, b, \eta^+\eta^-)
        = - \frac{b}{(1 - a (1 - \eta^+\eta^-))^2}
        ~~\text{and}~~
        \frac{\partial^2}{\partial (\eta^+\eta^-)^2} \,
        \varsigma(a, b, \eta^+\eta^-)
        =  \frac{2ab}{(1 - a(1-\eta^+\eta^-))^3}.
    \end{equation*}
    As $a \geq 0$ and $b \geq 0$ for all $(a,b) \in \set{A}$, the first and second derivatives are always nonpositive and nonnegative, respectively, which shows that $\varsigma$ is convex and nonincreasing in~$\eta^+\eta^-$. Since the pointwise maximum of convex functions is a convex function~\citep[p.~80]{SB04}, the asymptotic slope~$m$ is thus also convex and nonincreasing in~$\eta^+\eta^-$.
\end{proof}

\begin{proof}[Proof of Lemma~\ref{lem:lin}]
    If~$\dot y^\star = 0$, then we have
    \begin{equation*}
        \dot y(g(x^r), x^r) = 0 \iff \frac{\dot y(g(x^r), x^r)}{x^r} = \eta^+ \frac{g(x^r)}{x^r} - \eta_d \varphi\left( - \frac{g(x^r)}{x^r} \right) = 0 \iff g(x^r) = m x^r
    \end{equation*}
    for all~$x^r > 0$, where the second equivalence holds because~$m$ is the unique solution to $\eta^+ m - \eta_d \varphi(-m) = 0$, which is equivalent to $m = (1 - \eta^+\eta^-)\varphi(m)$, as revealed by the proof of Proposition~\ref{prop:if}. For $x^r = 0$, we have~$g(x^r) = 0 = m x^r$ because $g$ is continuous on~$\mathbb{R}_+$. Thus, $g(x^r) = m x^r$.
\end{proof}

\begin{proof}[Proof of Theorem~\ref{th:dr}] The robust constraints are replaced by their deterministic counterparts from Proposition~\ref{prop:cr}. The constraint on the expected terminal state-of-charge is enforced implicitly by expressing the decision variable~$x^b$ as the function~$g(x^r)$, characterized in Proposition~\ref{prop:if}. 
\end{proof}

\begin{proof}[Proof of Lemma~\ref{lem:cfs}]
    The feasible set~$\set{X}$ is convex if the function~$q(x^r) = g(x^r) - \ell(x^r)$ is monotonic, in which case there exists at most one intersection between $g$ and $\ell$. In the following, we show that Assumption~\ref{ass:ncnd_mad} implies that~$q$ is strictly decreasing. The slope of~$q$ is maximal when the slope of~$g$ is maximal and the slope of~$\ell$ is minimal. The maximal slope of~$g$ is~$m$, while the minimal slope of~$\ell$ is~$\frac{\gamma}{T}$. The function~$q$ is strictly decreasing if its maximal slope is strictly negative, which is the case if~$m < \frac{\gamma}{T}$. As~$\varphi(m) = \varphi(-m) + m$ by Lemma~\ref{lem:sym}, we have indeed
    \begin{equation*}
        m = (1 - \eta^+ \eta^-) \varphi(m)
        = \left( \frac{1}{\eta^+\eta^-} - 1 \right) \varphi(-m)
        \leq \left( \frac{1}{\eta^+\eta^-} - 1 \right) \varphi(0)
        < \frac{1}{2}\left( \frac{1}{\eta^+\eta^-} - 1 \right) \frac{\gamma}{T}
        \leq \frac{\gamma}{T}.
    \end{equation*}
    The first inequality holds because~$\varphi$ is nondecreasing and~$m \geq 0$. The strict inequality holds because~$\varphi(0) < \frac{\gamma}{2T}$ by Remark~\ref{rmk:mad}. The last inequality holds because~$\eta^+\eta^- \geq \frac{1}{3}$ by Assumption~\ref{ass:ncnd_mad}. 
\end{proof}

\begin{proof}[Proof of Theorem~\ref{th:x_ast}]
    Problem~\eqref{pb:P''} minimizes the net cost~$T(c^b g(x^r) - c^r x^r)$ over~$\set{X} = [0, \bar x^r]$. Since $g$ is convex by Proposition~\ref{prop:if}, its derivative~$g'$ is nondecreasing (where it exists). The net marginal cost $T(c^b g'(x^r) - c^r)$ is thus nondecreasing.
    As $g$ is also proper and closed, Theorem~24.1 by~\cite{RR70} implies that the right and left derivative functions $g'_+$ and $g'_-$ are nondecreasing and that $g'_+$ is a pointwise upper bound on~$g'_-$. If~$\bar x^r = 0$, the only feasible and therefore optimal solution is~$x^r_\ast = 0$. If $x^r > 0$, we distinguish three cases based on the values of $g'_+(0)$ and 
    $g'_-(\bar x^r)$.
    
    If~$g'_+(0) \geq \frac{c^r}{c^b}$, the marginal profit is nonpositive for all feasible values of~$x^r$, and $x^r_\ast = 0$ is the smallest optimal solution. 
    Conversely, if $g'_-(\bar x^r ) < \frac{c^r}{c^b}$, then the marginal profit is strictly positive for all feasible values of~$x^r$ and $x^r_\ast = \bar x^r$ is the only optimal solution. Finally, if~$g'(0)_+ < \frac{c^r}{c^b}$ and $g'_-(\bar x^r ) \geq \frac{c^r}{c^b}$, then the set~$\set{X_\star}$ of roots to the net marginal cost in~$(0,\bar x^r]$ is nonempty and compact because $g$ is convex and continuous. Any root would be an optimal solution and, in particular, the smallest root $x^r_\ast = \min \set{X}_\star$ is an optimal solution. This solution exists because~$\set{X_\star}$ is compact.    
\end{proof}

\begin{proof}[Proof of Lemma~\ref{lem:phi}]
As~$\tilde \xi$ is supported on~$[-1,1]$ we have \mbox{$\ubar \varphi (\xi) = \varphi(\xi) = \bar \varphi(\xi)$} for all $\vert \xi \vert \geq 1$. For any~$\vert \xi \vert < 1$, $\ubar \varphi(\xi) \leq \varphi(\xi)$ as $\ubar \varphi$ is the piecewise maximum of three affine functions that are tangent to the convex function~$\varphi$ at~$\xi = -1$, $ = 0$, and~$ = 1$. Conversely, $\bar \varphi(\xi) \geq \varphi(\xi)$ since $\bar \varphi$ is the piecewise maximum of two linear interpolations of~$\varphi$ from~$\xi = -1$ to $\xi = 0$ and from~$\xi = 0$ to~$\xi = 1$.
\end{proof}

\begin{proof}[Proof of Lemma~\ref{lem:g_bds}]
    For any~$x^r \geq 0$, the constraint~$\dot y(x^b, x^r) = \dot y^\star$ implies that $\ubar g(x^r)$, $g(x^r)$, and $\bar g(x^r)$ are the unique roots of the functions~$s_{\ubar \varphi}$, $s_\varphi$, and $s_{\bar \varphi}$, respectively, where~$s_\varphi$ is defined as
    \begin{equation*}
        s_\varphi(x^b) = \dot y^\star + \eta_d \varphi\left(-\frac{x^b}{x^r}\right)x^r - \eta^+ x^b
    \end{equation*}
    and the functions~$s_{\ubar \varphi}$ and~$s_{\bar \varphi}$ are defined similarly. By definition, we have~$0 = s_{\ubar \varphi}(\ubar g(x^r)) \leq s_\varphi( \ubar g(x^r))$, where the inequality holds because $\ubar \varphi(\ubar g(x^r)) \leq \varphi(\ubar g(x^r))$, see Lemma~\ref{lem:phi}, and $\eta_d x^r \geq 0$. Since $s_{\ubar \varphi}$ and $s_\varphi$ are strictly decreasing by Proposition~\ref{prop:esoc}, we must have $g(x^r) \geq \ubar g(x^r)$. A similar reasoning shows that $\bar g(x^r) \geq g(x^r)$. As $\ubar g$, $g$, and $\bar g$ are convex, the inequality $\ubar g(x^r) \leq g(x^r) \leq \bar g(x^r)$ implies the inequality on the asymptotic slopes $\underline{m} \leq m \leq \overline{m}$.
\end{proof}

\begin{proof}[Proof of Proposition~\ref{prop:xr_lin}]
    If $y_0 = y^\star$, then $\dot y^\star = \frac{y^\star - y_0}{T} = 0$ and so $g(x^r) = m x^r$ by Lemma~\ref{lem:lin}. Hence, the constraints in Problem~\eqref{pb:P''} obey the following equivalences.
    \begin{align*}
        x^r + g(x^r) \leq \bar y ^+ & 
        \quad \iff \quad 
        x^r \leq \frac{\bar y^+}{1+m} \\
        x^r - g(x^r) \leq \bar y^- &
        \quad \iff \quad
        x^r \leq \frac{\bar y^-}{1-m} \\
        x^r + \max\left\{ \frac{T}{\gamma}g(x^r), g(x^r) \right\} \leq \frac{\bar y - y_0}{\eta^+ \gamma} &
        \quad \iff \quad
        x^r \leq \frac{\bar y - y_0}{\eta^+(\gamma + mT)}\\
        x^r - \min\left\{ \frac{T}{\gamma} g(x^r), g(x^r) \right\} \leq \frac{\eta^- y_0}{\gamma} & 
        \quad \iff \quad
        x^r \leq \frac{\eta^- y_0}{\gamma(1-m)}
    \end{align*}
    The last two equivalences hold because~$\frac{T}{\gamma} \geq 1$ and $g(x^r) = m x^r \geq 0$ since any feasible~$x^r$ must be nonnegative and since $m$ is nonnegative by Proposition~\ref{prop:if}.
\end{proof}

\begin{proof}[Proof of Theorem~\ref{th:as}]
    If $y_0 = y^\star$, then $\dot y^\star = \frac{y^\star - y_0}{T} = 0$ and so $g(x^r) = m x^r$ by Lemma~\ref{lem:lin}. Hence, $g'_+(0) = g'_-(\bar x^r) = m$. Theorem~\ref{th:x_ast} implies that $x^r_\ast = 0$ if $m \geq \frac{c^r}{c^b}$ and $= \bar x^r$ otherwise.
\end{proof}
\newpage
\section{Supplementary Material}

As supplementary material, we first provide additional descriptions and analysis of the case study presented in Section~\ref{sec:appl} of the main paper. Next, we present a model that captures elastic prices. Finally, we provide all proofs that are not included in Appendix~\ref{sec:ch3_apx_proofs} of the main paper.

\subsection{Detailed Model Parameters}\label{sm:params}

\subsubsection{Frequency Regulation}
The technical term for the type of frequency regulation we consider is \emph{frequency containment reserves}~(FCR). France participates in a common European market for frequency containment reserves with a daily planning horizon.\footnote{\url{https://entsoe.eu/network_codes/eb/fcr}} We thus set~$T = 24$~hours. In its regulation on frequency containment reserves the European Commission specifies that the \emph{``minimum activation period to be ensured by FCR providers \textnormal{[is not to be]} greater than 30 or smaller than 15~minutes''} and that storage operators \emph{``shall ensure the recovery of \textnormal{[their]} energy reservoirs as soon as possible, within 2~hours after the end of the alert state''} \cite[art.~156(10, 13)]{EU17}, where an \emph{activation period} designates a period of consecutive extreme frequency deviations~$\delta(t) \in \{-1,1\}$. The total activation period~$\gamma$ to be ensured over a period of 24~hours is thus between 2.75~hours and 5~hours. The uncertainty set~$\set{D}$ contains all frequency deviation signals that correspond to a given total activation period~$\gamma$. Some of these signals may exhibit activation periods that are longer than the minimum activation period prescribed by the European Commission. The uncertainty set~$\set{D}$ is therefore a conservative approximation of the regulation by the European Commission. The strength of the delivery guarantee required by the European Commission can nevertheless be measured by the \emph{activation ratio}~$\gamma/T$, that is, the fraction of time during which a storage operator must be able to provide all the regulation power she promised. When studying the profits a storage operator can reap from frequency regulation over the lifetime of an electricity storage device, we will vary the length of the planning horizon~$T$ while keeping the activation ratio constant. In line with the regulation by the European Commission, we will consider activation ratios of $0.1$ and $0.2$.

\subsubsection{Electricity Prices}
\begin{figure}[t!]
	\centering
    \includegraphics[width=1\linewidth]{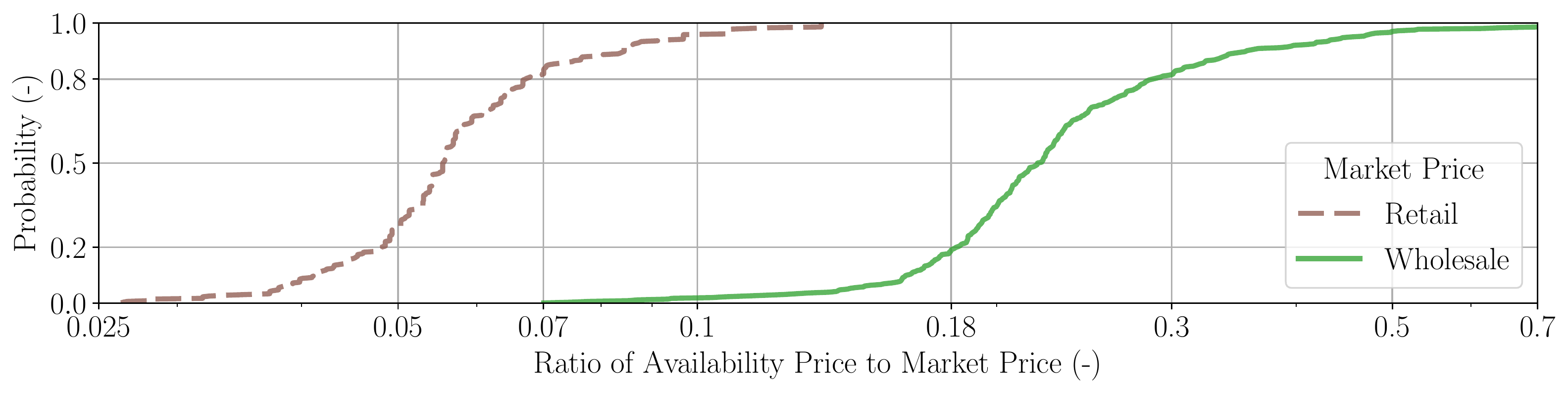}
    \caption{Distribution of daily average electricity prices in 2019.}
    \label{fig:cr_cb}
\end{figure}

In the year 2019, the expected average price of regulation power was the same as the expected average availability price because $\mathbb{E} \frac{1}{T}\int_\set{T} \tilde \delta(t) \tilde p^d(t) \, \mathrm{d}t$ vanished. In fact, the average value of $\frac{1}{T}\int_\set{T} \delta(t) p^d(t) \, \mathrm{d}t$ over all days was~$-8.77\cdot 10^{-5}$. The minimum availability, wholesale market, and retail market prices in any half-hour interval were~$0.41\frac{\text{cts}}{\text{kW}\cdot \text{h}}$, $-2.49\frac{\text{cts}}{\text{kWh}}$, and $14.5\frac{\text{cts}}{\text{kWh}}$, respectively. Averaged over each day, the minimum daily average availability, wholesale market, and retail market prices were~$0.41\frac{\text{cts}}{\text{kW}\cdot\text{h}}$, $0.37\frac{\text{cts}}{\text{kWh}}$, and $14.5\frac{\text{cts}}{\text{kWh}}$, respectively, which are all strictly positive. The ratio of daily average availability prices to daily average market prices of electricity was between~$0.070$ and~$2.168$ with an average of~$0.251$ for wholesale market prices, and between~$0.026$ and~$0.133$ with an average of~$0.059$ for retail market prices. Figure~\ref{fig:cr_cb} shows the empirical cumulative distribution function of this ratio for wholesale and retail market prices. When studying the profits from frequency regulation over the planning horizon~$\set{T}$ and over the lifetimes of electricity storage devices, we will set the ratio of the expected average price of regulation power~$c^r$ to the expected average market price of electricity~$c^b$ to~$\frac{c^r}{c^b} = 0.251$ for wholesale market prices and to~$\frac{c^r}{c^b} = 0.059$ for retail market prices.

The desired charging rate~$\dot y^\star$ for meeting the terminal state-of-charge target influences the quantity and the price of regulation power that a storage operator can offer. Ideally, the storage operator would be able to meet the terminal state-of-charge target~$y^\star$ not just in expectation but exactly. In the following, we assume that this is the case and that the terminal state-of-charge target stays constant from one planning horizon to another. The desired charging rate vanishes therefore during any given planning horizon. In this case, the implicit function $g$ is the linear function $g(x^r) = mx^r$ by Proposition~\ref{prop:if}. The marginal cost of providing frequency regulation is thus $Tmc^b$. If the desired charging rate was nonzero, then the marginal cost would be lower as the slope of~$g$ would converge to~$m$ only asymptotically. We thus overestimate the marginal cost of providing frequency regulation when assuming that $\dot y^\star = 0$. Regardless of the particular value of~$\dot y^\star$, it is therefore always profitable for the storage operator to sell as much regulation power as possible if the marginal revenue $Tc^r$ of providing frequency regulation is higher than the marginal cost $Tmc^b$ given~$\dot y^\star = 0$, \ie, if $\frac{c^r}{c^b} > m$.

\subsubsection{Frequency Deviation Distribution}\label{sm:delta}
Over the years 2017 to 2019, the cumulative distribution function corresponding to~$\mathbb{P}_\xi$ can be approximated by a symmetric logistic function~$F$ with a maximum error of~$0.018$ with respect to the empirical cumulative distribution function constructed from about $9.5$~million frequency recordings with a 10~second resolution. This justifies Assumption~\ref{ass:phi}. Based on the logistic approximation, the cumulative distribution function and the characteristic function are given by $F(\xi) = \frac{1}{1+\exp{(-\theta \xi)}}$ and $\varphi(\xi) = \frac{\ln({1 + \exp{(\theta \xi)}})}{\theta}$, respectively, where $\theta = \frac{2\ln{(2)}}{\Delta}$ and $\Delta = 2\varphi(0) = 0.0816$, which satisfies the condition $\Delta \leq \frac{\gamma}{T}$ in Assumption~\ref{ass:ncnd_mad} as~$\frac{\gamma}{T} \geq 0.1$. The coefficient~$\theta$ was chosen such that frequency deviations have the same mean absolute deviation~$\Delta$ under the logistic distribution as under the empirical distribution. In principle, $F$ should be a truncated logistic function as the support of~$\mathbb{P}_\xi$, $\Xi = [-1,1]$, is bounded. The truncation error, however, is just $2.5\cdot 10^{-8}$, which we deem to be negligible. We find that the mean absolute deviation of the frequency recordings is smaller than~$0.1$ on $96.7\%$ of all days. If the uncertainty set~$\set{D}$ is parametrized by~$\gamma = 0.1T$, the empirical frequency deviation signals fall thus outside of~$\set{D}$ on~$3.3\%$ of all days. On these days, the storage operator may stop delivering regulation power once the total absolute deviation~$\int_0^t \vert \delta(t') \vert \, \mathrm{d}t'$ of the frequency deviation signal exceeds~$\gamma$. In principle, the frequency deviation distribution should thus be estimated based on past frequency deviation signals with mean absolute deviation capped at the activation ratio~$\gamma/T$. The distribution estimated directly on past frequency deviation signals will differ most from the distribution estimated on capped past frequency deviation signals if the activation ratio is equal to~$0.1$ rather than~$0.2$. In this case, the maximum difference between the two distributions is less than~$0.001$, which we consider negligible compared to the maximum difference of~$0.018$ between the empirical cumulative distribution function and the logistic function~$F$.

\begin{figure}[t!]
    \centering
    \includegraphics[width=1\linewidth]{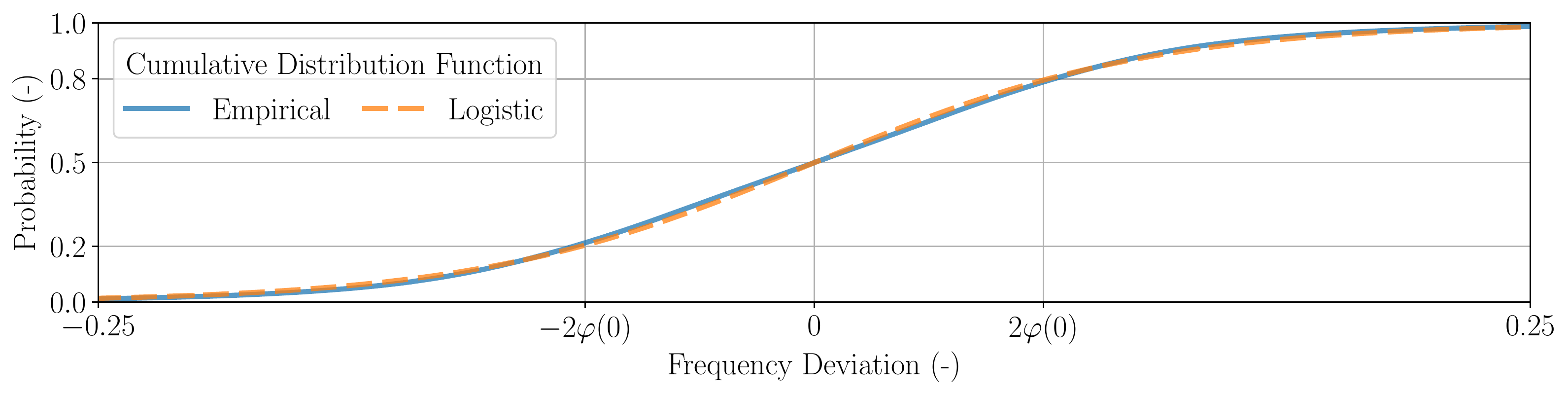}
    \caption{Distribution of frequency deviations in 2017--19.}
    \label{fig:F_sig}
\end{figure}

Figure~\ref{fig:F_sig} shows the empirical cumulative distribution function and its logistic approximation.

\subsection{Detailed Results}

\subsubsection{Impact of Charging and Discharging Losses on the Amount of Regulation Power}\label{sm:xr}

Charging and discharging losses may  impact the amount of regulation power~$\bar x^r$ a storage device can provide. We have seen in Section~\ref{sec:as} that, if~$\dot y^\star = 0$, then the storage device may be either energy-constrained or power-constrained.

If the storage device is power-constrained, then the storage operator can account for the roundtrip efficiency by dimensioning the discharging capacity~$\bar y^-$ as a fraction~$\frac{1-m}{1+m}$ of the charging capacity~$\bar y^+$ without restricting~$\bar x^r$. The fraction is~$1$ at a roundtrip efficiency of~$1$, and decreases to~$0.92$ as the roundtrip efficiency decreases to~$0.35$.

If the storage device is energy-constrained, then charging losses increase the amount of energy that the storage device can consume from the grid, while discharging losses decrease the amount of energy that the storage device can deliver to the grid. The storage operator can offer most regulation power $\bar x^r_\star$ to the grid operator if the initial state-of-charge is such that she can consume as much energy from the grid as she can provide to the grid. For roundtrip efficiencies in~$(0,1)$, the optimal initial state-of-charge $y_0^\star$ is nondecreasing in the activation ratio. Given activation ratios between~$0.1$ and~$0.2$, $\frac{y^\star_0}{\bar y}$ is between $0.52$ and $0.53$ for lithium-ion batteries with a roundtrip efficiency of~$0.85$, between $0.57$ and $0.60$ for redox flow batteries with a roundtrip efficiency of~$0.60$, and between $0.66$ and~$0.69$ for hydrogen storage with a roundtrip efficiency of~$0.35$.

\begin{figure}[t!]
\centering
\begin{subfigure}{0.49\textwidth}
    \includegraphics[width=\textwidth]{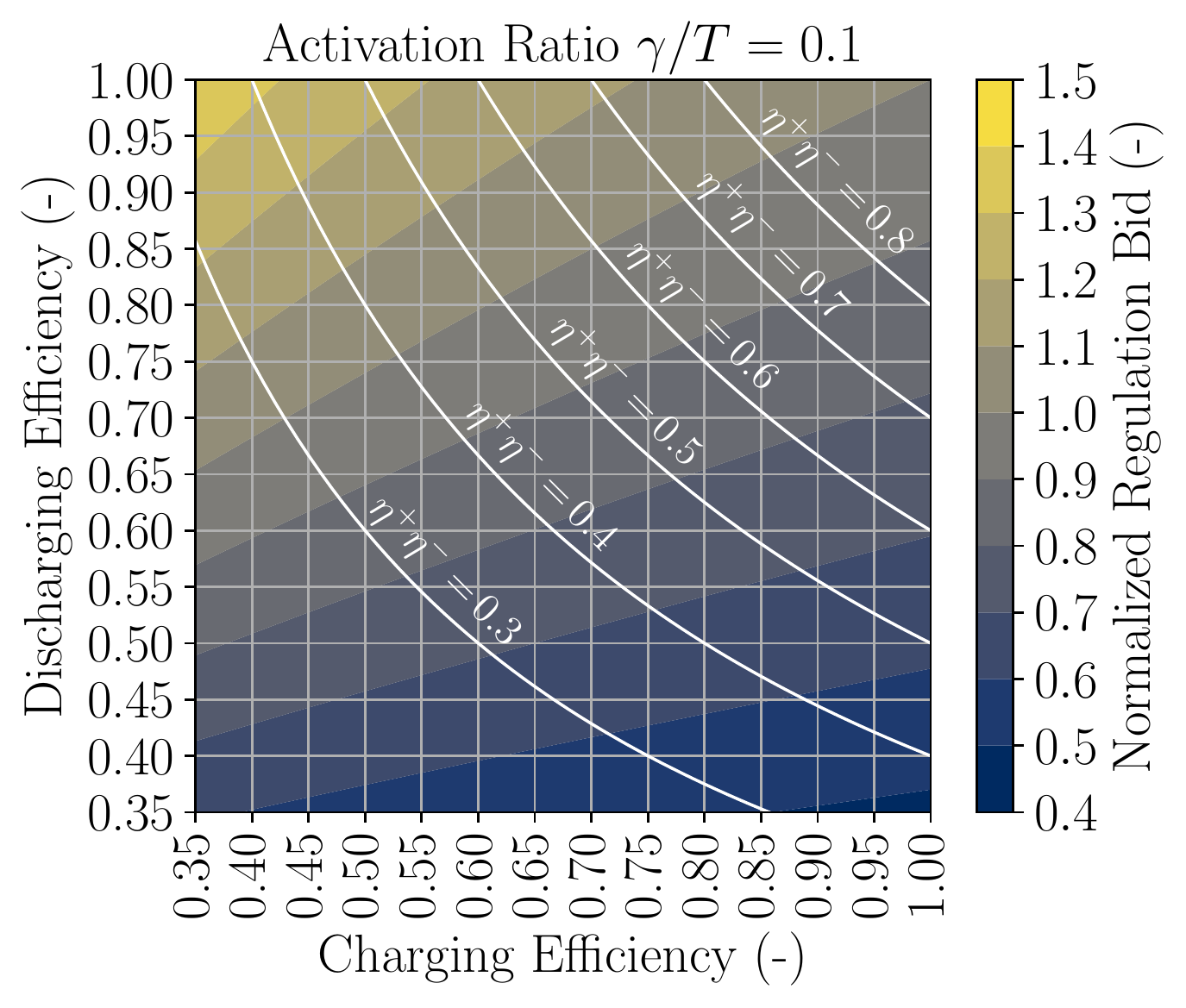}
\end{subfigure}
\hfill
\begin{subfigure}{0.49\textwidth}
    \includegraphics[width=\textwidth]{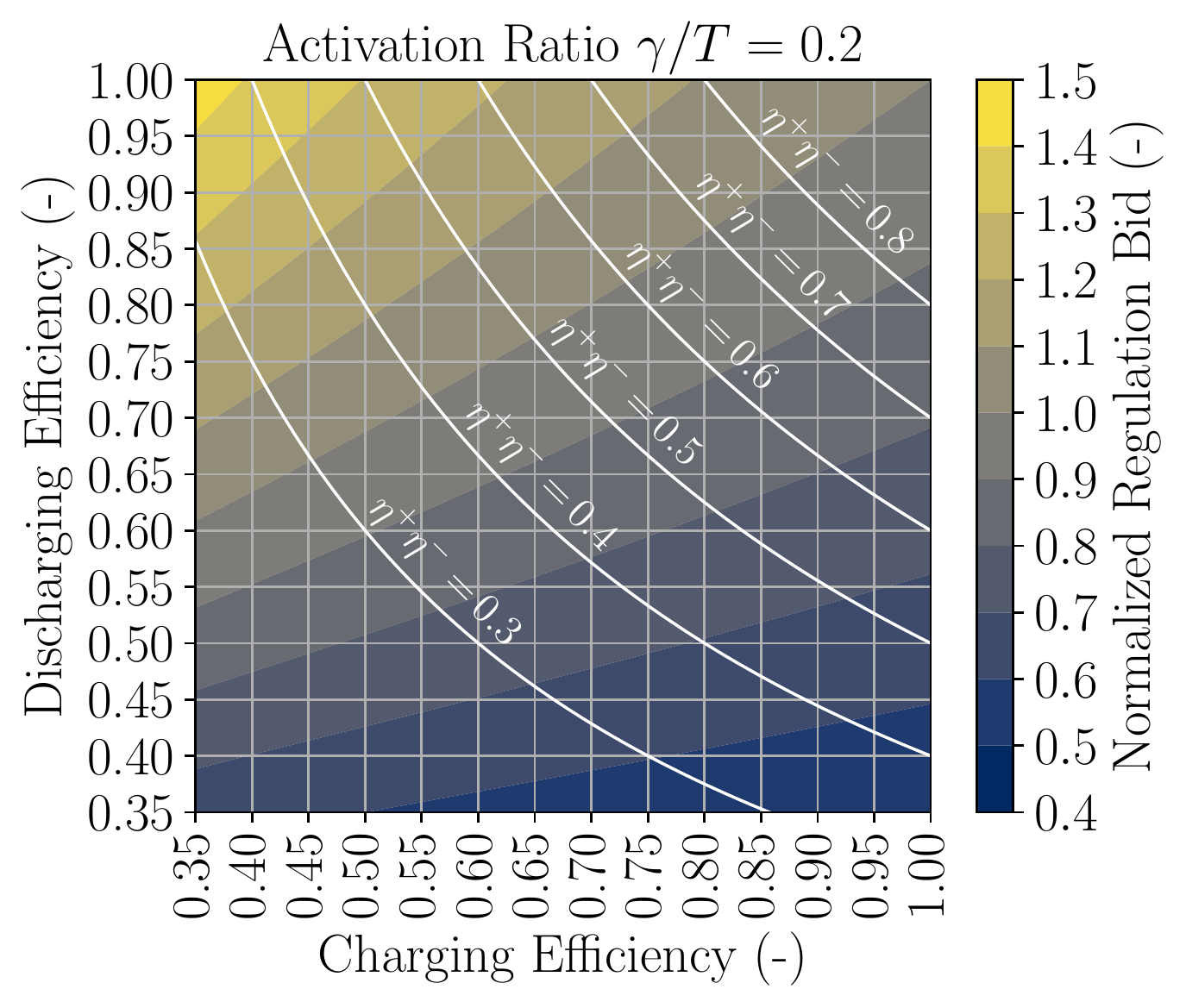}
\end{subfigure}
\caption{Maximum regulation power normalized by~$\frac{\bar y}{2\gamma}$.}
\label{fig:xr_max}
\end{figure}

Figure~\ref{fig:xr_max} shows that the normalized regulation power~$\bar x^r_\star / \frac{\bar y}{2\gamma}$ is indeed nonincreasing in the charging efficiency~$\eta^+$ and nondecreasing in the discharging efficiency~$\eta^-$. The decrease in~$\eta^+$ is more pronounced and the increase in~$\eta^-$ is less pronounced if the activation ratio is~$0.2$ rather than~$0.1$. Starting from a roundtrip efficiency of~$1$ and an activation ratio of~$0.2$, the normalized regulation power increases from~$1$ to~$1.45$ as the charging efficiency decreases from~$1$ to~$0.35$, and decreases from~$1$ to~$0.51$ as the discharging efficiency decreases from~$1$ to~$0.35$. At an activation ratio of~$0.1$, the normalized regulation power increases to only~$1.37$ as the charging efficiency decreases to~$0.35$, and decreases slightly further to~$0.48$ as the discharging efficiency decreases to~$0.35$. Although the normalized regulation power is lower, the absolute regulation power is considerably higher at an activation ratio of~$0.1$ rather than~$0.2$, because the normalization constant $\frac{\bar y}{2\gamma}$ is inversely proportional to~$\gamma$. 

\subsubsection{Operating Profits}\label{sm:op}

\begin{figure}[t!]
\centering
\begin{subfigure}{0.49\textwidth}
    \includegraphics[width=\textwidth]{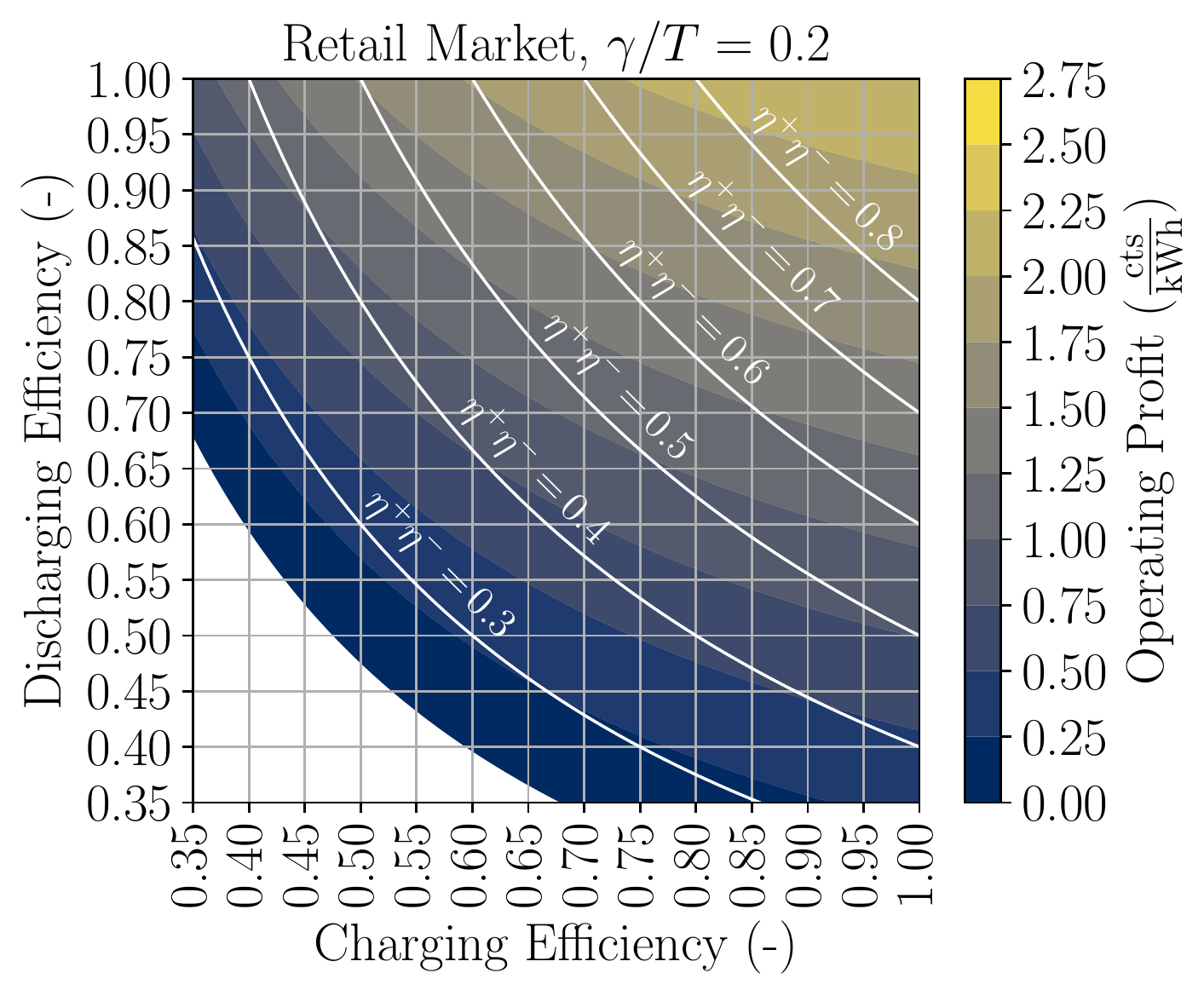}
\end{subfigure}
\hfill
\begin{subfigure}{0.49\textwidth}
    \includegraphics[width=\textwidth]{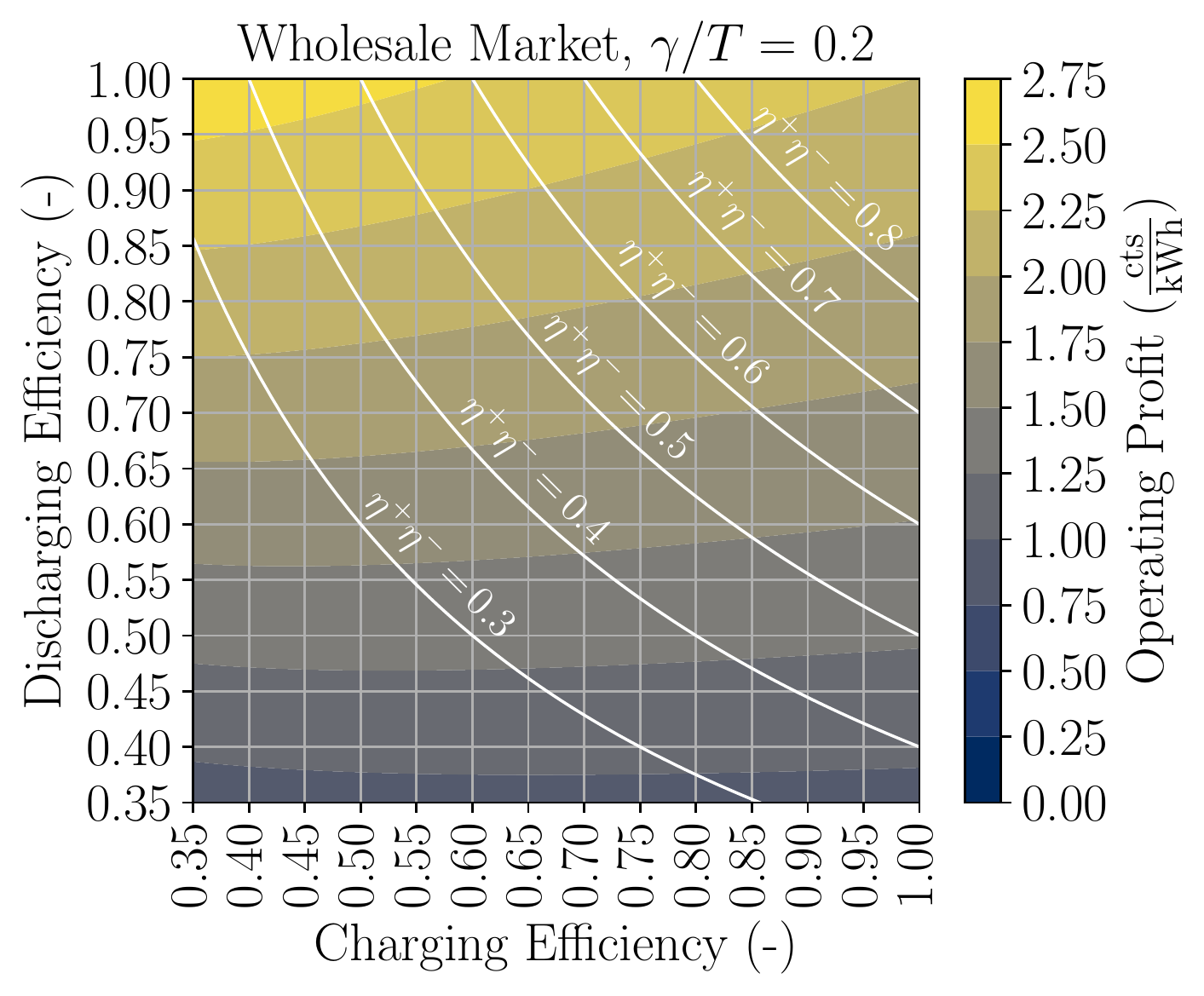}
\end{subfigure}
\begin{subfigure}{0.49\textwidth}
    \includegraphics[width=\textwidth]{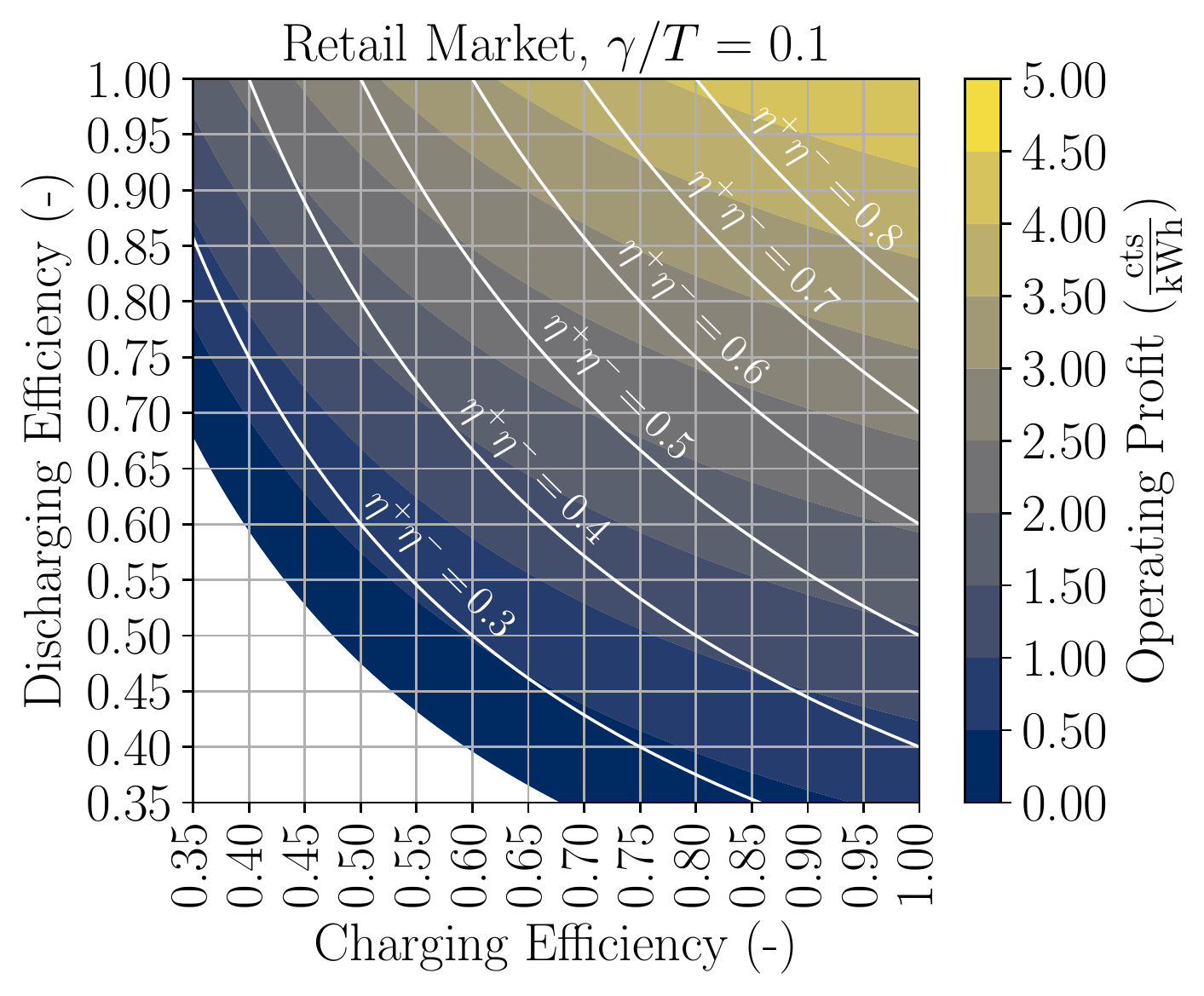}
\end{subfigure}
\hfill
\begin{subfigure}{0.49\textwidth}
    \includegraphics[width=\textwidth]{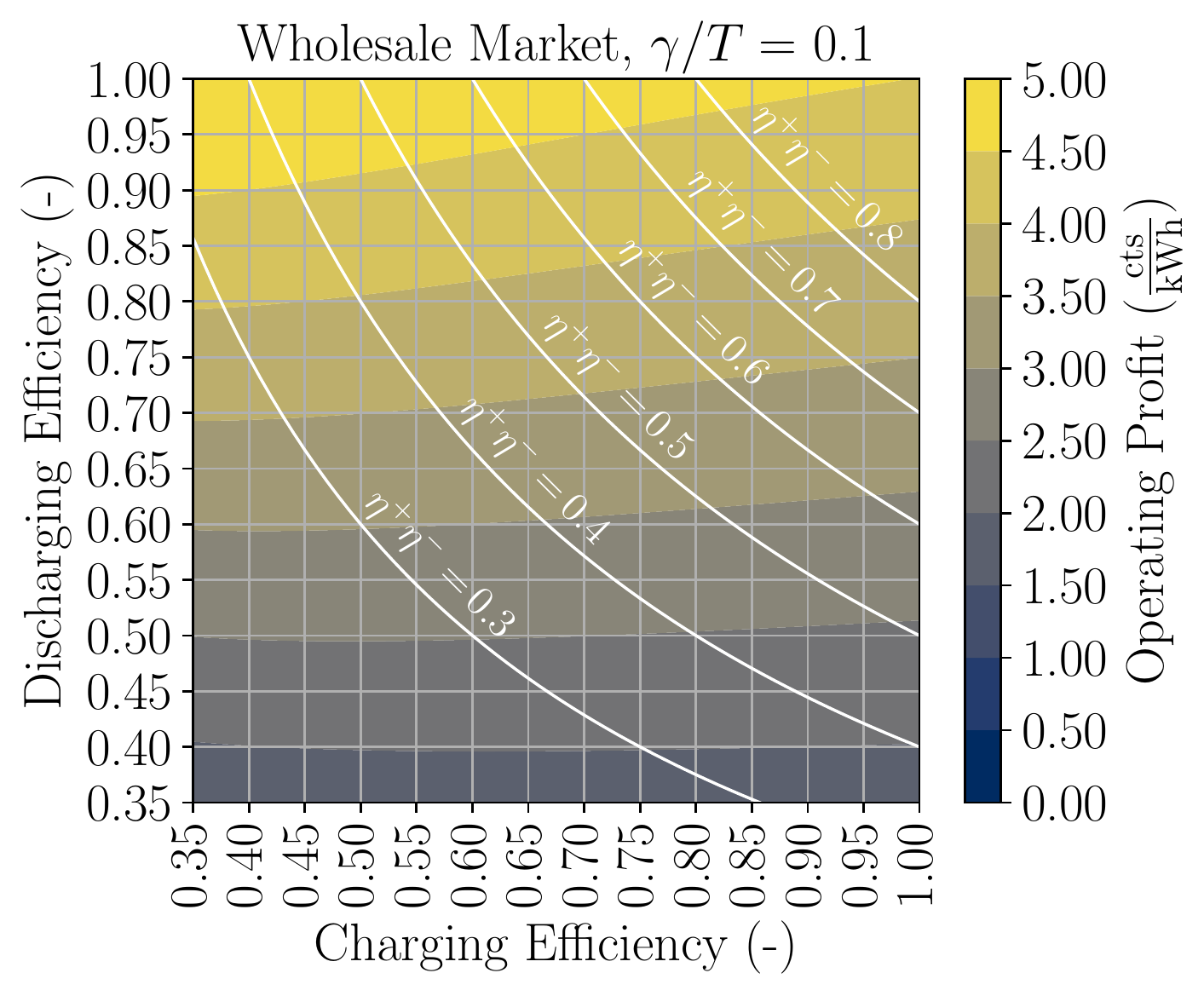}
\end{subfigure}
\caption{Operating profit per kWh of storage capacity at wholesale and retail prices.}
\label{fig:pi_T}
\end{figure}

We established in Section~\ref{sec:m_xr} that the profit per unit of regulation power increases with the roundtrip efficiency, while the maximum amount of regulation power increases with the discharging efficiency but decreases with the charging efficiency. The operating profit thus increases with the discharging efficiency, but it is unclear how it depends on the charging efficiency. Intuitively, the higher the market price of electricity, the more pronounced the increase of the profit per unit of regulation power in the charging efficiency. Similarly, the higher the activation ratio, the more pronounced the decrease of the maximum amount of regulation power in the charging efficiency. Figure~\ref{fig:pi_T} shows that, for activation ratios of~$0.2$ and~$0.1$, the operating profit increases with the charging efficiency at retail electricity prices but not at wholesale electricity prices. In practice, however, the combined effect of charging and discharging losses is usually a decrease in the operating profit, even at wholesale prices. As examples, we consider again lithium-ion batteries, vehicle-to-grid, and hydrogen storage, with the same charging and discharging efficiencies as in Section~\ref{sec:m_xr}. At an activation ratio of~$0.2$, the operating profit is~$2.25\frac{\text{cts}}{\text{kWh}}$ in the absence of charging and discharging losses, regardless of the market price of electricity. For lithium-ion batteries, the operating profit reduces to~$2.15\frac{\text{cts}}{\text{kWh}}$ at wholesale prices and to~$1.96\frac{\text{cts}}{\text{kWh}}$ at retail prices. For vehicle-to-grid, the operating profit decreases further to~$1.92\frac{\text{cts}}{\text{kWh}}$ at wholesale prices and to~$1.53\frac{\text{cts}}{\text{kWh}}$ at retail prices. For hydrogen, finally, the operating profit falls to~$1.49\frac{\text{cts}}{\text{kWh}}$ at wholesale prices and to~$0.81\frac{\text{cts}}{\text{kWh}}$ at retail prices. If the activation ratio halves from~$0.2$ to~$0.1$, the operating profits roughly double.

\subsection{Model with elastic prices}\label{sec:si_elastic_prices}

If we allow for elastic prices that depend on the amount of power traded on the markets, then the expected cost over the planning horizon becomes
\begin{align*}
        & \E \int_\T \tilde p^b(x^b, t) x^b - \left( \tilde p^a(x^r, t) - \tilde \delta(t) \tilde p^d(x^b, t) \right) x^r \, \mathrm{d}t 
        = 
        \E \int_\T \tilde p^b(x^b, t) x^b - \tilde p^a(x^r, t) x^r \, \mathrm{d}t \\
        = &
        T \left( c^b_0  x^b + c^b_d (x^b)^2 - c^a_0 x^r + c^a_d (x^r)^2\right)
        = 
        T \left( c^b_0  g(x^r) + c^b_d g(x^r)^2 - c^a_0 x^r + c^a_d (x^r)^2\right),
\end{align*}
where we assumed that delivery prices are independent of frequency deviations, so the first equality holds. Following \cite{barbry2019robust}, we also assume that market prices are affine functions of the amount of power traded on markets, so the second equality holds. The third equality follows from $x^b = g(x^r)$. To model decreasing marginal prices, the cost function must be convex in $x^b$ and $x^r$, hence $c^b_d \geq 0$ and $c^a_d \geq 0$. In addition, for the prices to be positive we must have $c^b_0 > 0$ and $c^a_0 > 0$.

\begin{Prop}
    If $g(0) \geq -\frac{c^b_0}{2c^b_d}$, the function $T \left( c^b_0  g(x^r) + c^b_d g(x^r)^2 - c^a_0 x^r + c^a_d (x^r)^2\right)$
    is convex.
\end{Prop}

\begin{proof}
    Let $h(x^b) = c^b_0  x^b + c^b_d (x^b)^2$ and $f = h \circ g$. Then, the cost function can be written as $ T \left( f(x^r) - c^a_0 x^r + c^a_d (x^r)^2\right)$, which is convex if $f$ is convex because nonnegative weighted sums preserve convexity~\citep[p.~79]{SB04}. By the composition rules in \citeauthor{SB04}~(2004, p.~84), $f$ is convex if $h$ is convex and nondecreasing, and $g$ is convex. We know that $h$ is convex as $c^b_d \geq 0$ and, from Proposition~\ref{prop:if}, that $g$ is convex; $h$ is nondecreasing if
    \begin{equation*}
        h'(g(x^r)) = c^b_0 + 2 c^b_d g(x^r) \geq 0~~\forall x^r \in \mathbb{R}_+ \iff g(0) \geq - \frac{c^b_0}{2c^b_d},
    \end{equation*}
    where the equivalence holds because $g$ is nondecreasing. If $g(0) \geq -\frac{c^b_0}{2c^b_d}$, $f$ is thus indeed convex. 
\end{proof}

As $\ell(0) \leq g(0)$ for any feasible solution, we know that the cost function is convex if $-\frac{c^b_0}{2c^b_d} \leq \ell(0) = -\min\left\{ \eta^- \frac{\bar y}{T}, \bar y^- \right\}$. Typically, this condition is fulfilled as $c^b_d$ will be small relative to~$c^b_0$. 

\begin{figure}
    \centering
    \includegraphics[width = \linewidth]{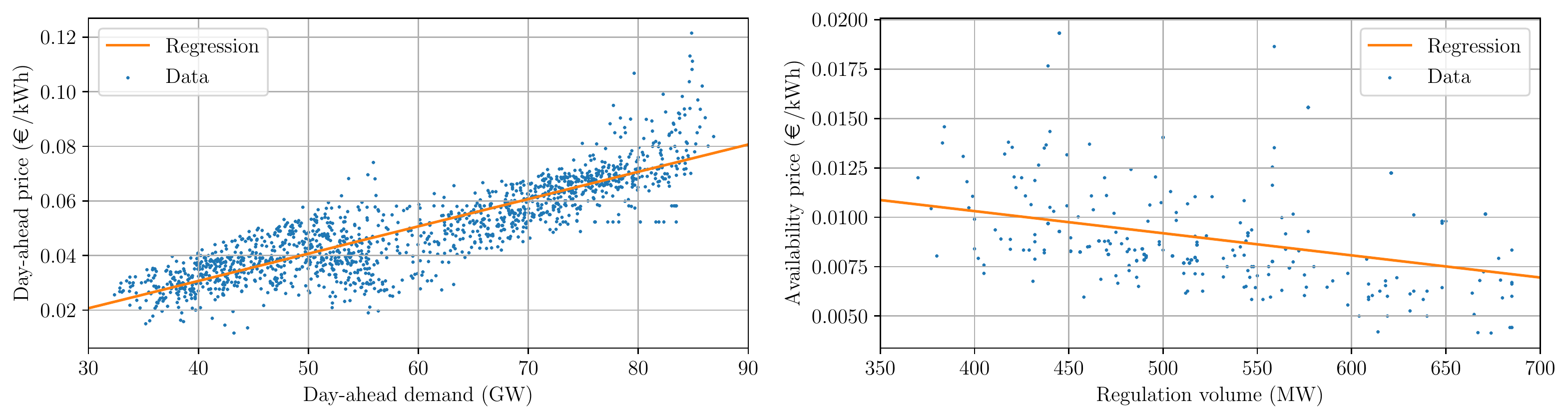}
    \caption{Day-ahead data with hourly resolution for January and July 2019, primary frequency regulation data with daily resolution for the full year 2919.}
    \label{fig:elastic_prices}
\end{figure}

\begin{Ex}~\label{ex:elastic_prices_cb}
    Consider a storage operator that trades electricity on the French day-ahead market. In the year 2019, the mean spot market price was $c^b_0 = 3.9\frac{\text{cts}}{\text{kWh}}$ and a least-squares linear regression yields a price elasticity of $c^b_d \approx 10^{-7}\frac{\text{cts}}{\text{kWh}\cdot\text{kW}}$ (see Figure~\ref{fig:elastic_prices}). Hence, $\frac{c^b_0}{2c^b_d} \approx 19.5$GW. In comparison, the total installed capacity of the European pumped hydro storage is about 50GW~\citep{EU20}. The condition $\frac{c^b_0}{2c^b_d} \geq \min\left\{ \eta^- \frac{\bar y}{T}, \bar y^- \right\}$ is thus very likely to hold for most storage operators.
\end{Ex}

In the remainder, we will assume that the inequaltiy holds strictly.
\begin{Ass}\label{ass:elasticity}
    We have that $\frac{c^b_0}{2c^b_d} > \min\left\{ \eta^- \frac{\bar y}{T}, \bar y^- \right\}$.
\end{Ass}

We will now examine the impact of elastic prices on the structure of the optimal solution to problem~\eqref{pb:P}. The marginal cost, where it exists, becomes
\begin{equation*}
    T \left( c^b_0 g'(x^r) + 2 c^b_d g'(x^r) g(x^r) - c^a_0 + 2 c^a_d x^r \right),
\end{equation*}
which is nondecreasing because the cost function is convex under Assumption~\ref{ass:elasticity}. Intuitively, the marginal cost should indeed be nondecreasing because we expect ($i$) that the more regulation power is sold, the higher the losses and the more power needs to be bought from the market to cover them, and ($ii$) that the more regulation power is sold, the lower the availability price and the higher the price of buying electricity from the market. Expectation ($ii$) only applies when prices are elastic, expectation ($i$) applies more generally. The set of stationary points becomes
\begin{equation}\label{eq:Xstar_elastic}
    \set{X}_\star = \left\{
        x^r \in [0, \bar x^r]: \frac{c^a_0 - 2 c^a_d x^r}{c^b_0 + 2 c^b_d g(x^r)} \in \partial g(x^r)
    \right\},
\end{equation}
which gives rise to the following modified version of Theorem~\ref{th:x_ast}.
\begin{Th}\label{th:x_ast_elastic}
    For $\bar x^r \geq 0$, the smallest optimal solution to problem~\emph{(P)} with elastic prices is
    \begin{equation*}
        x^r_\ast = \begin{cases}
         0 & \text{if } \bar x^r = 0 \text{ or } \bar x^r > 0 \text{ and } g'_+(0) \geq \frac{c^a_0}{c^b_0 + 2 c^b_d g(0)}, \\
        \bar x^r & \text{if } \bar x^r > 0 \text{ and } g'_-(\bar x^r) < \frac{c^a_0 - 2 c^a_d \bar x^r}{c^b_0 + 2 c^b_d g(\bar x^r)}, \\
        \min \set{X}_\star & \text{otherwise}.
        \end{cases}
    \end{equation*}
\end{Th}

\begin{proof}
    Substitute $x^r = 0$ and $x^r = \bar x^r$ into the condition for the stationary points.
\end{proof}

Remark~\ref{rmk:inf_X} applies to Theorem~\ref{th:x_ast_elastic} because the marginal cost is nondecreasing. We now focus on the case $y_0 = y^\star$ to derive a counterpart to Theorem~\ref{th:as}.

\begin{Th}\label{th:as_elastic}
If~$y_0 = y^\star$, then an optimal solution to problem~\emph{(P)} with elastic prices is
\begin{equation*}
    x^r_\ast = \begin{cases}
    0 & \text{if } m \geq \frac{c^a_0}{c^b_0}, \\
    \bar x^r & \text{if } \frac{c^a_0 - m c^b_0}{2( c^a_d + m^2 c^b_d )} \geq \bar x^r, \\
    \frac{c^a_0 - m c^b_0}{2( c^a_d + m^2 c^b_d )} & \text{otherwise},
    \end{cases}
\end{equation*}
under any frequency deviation distribution~$\mathbb{P}_\xi$. If~$\mathbb{P}_\xi = \ubar{\mathbb{P}}_\xi$, then $m = \underline{m}$.  If~$\mathbb{P}_\xi = \bar{\mathbb{P}}_\xi$, then $m = \overline{m}$.
\end{Th}

\begin{proof}
Substituting $g(x^r) = m x^r$ and $\partial g(x^r) = m$ into~\eqref{eq:Xstar_elastic} yields
    \begin{equation*}
        \set{X}_\star = \left\{
        x^r \in [0, \bar x^r]: x^r = \frac{c^a_0 - m c^b_0}{2( c^a_d + m^2 c^b_d )}
    \right\}. \qedhere
    \end{equation*}
\end{proof}

\begin{table}[b]
    \centering
    \begin{tabular}{crccc}
        Prices $(-)$ & Effective yearly profit (\EUR{}/kWh) & $T$ (h) & $\gamma/T$ $(-)$ & Investment costs  \\
        \midrule
        Inelastic & $1.600$ & $24$ & $0.1$ & Low \\
        Elastic   & $1.592$ & $24$ & $0.1$ & Low \\
        \midrule
        Inelastic & $50.745$ & $4$ & $0.1$ & Low \\
        Elastic   & $50.693$ & $4$ & $0.1$ & Low \\
        \bottomrule
    \end{tabular}
    \caption{Effective yearly profit per kWh of storage capacity for a 100kWh battery.}
    \label{tab:pi_y_elastic}
\end{table}

Under inelastic prices there were only two candidates for the optimal solution: offering no regulation power at all or as much regulation power as possible. Under elastic prices, there is an additional candidate: offering an intermediate amount of regulation power.

\begin{Ex}\label{ex:power}
    Consider the same storage operator as in Example~\ref{ex:elastic_prices_cb}. In the year 2019, the mean availability price was $c^a_0 = 0.9\frac{\text{cts}}{\text{kWh}}$ and a least-squares linear regression yields a price elasticity of $c^a_d \approx 10^{-6}\frac{\text{cts}}{\text{kWh}\cdot\text{kW}}$ (see Figure~\ref{fig:elastic_prices}). Assume $m = 0.03$, which is representative of storage devices with a roundtrip efficiency of~$0.50$ (see Figure~\ref{fig:m}). Then, $\frac{c^a_0 - m c^b_0}{2( c^a_d + m^2 c^b_d)} \approx 450$MW, which is close to the mean amount of frequency regulation, $520$MW, bought by the French grid operator in the year~2019. The third candidate solution is thus optimal only if the storage operator supplies almost the entire need for regulation power of the French grid operator. 
\end{Ex}

To judge how elastic prices may impact our conclusions about the impact of a reduced planning horizon (see Figure~\ref{fig:pi_y}), we now compare the effective yearly profit of investing in a 100kWh battery if the planning horizon was reduced from 24 hours to 4 hours under elastic and inelastic prices. Table~\ref{tab:pi_y_elastic} shows that the profits under elastic prices are about $0.1$\% to $0.5$\% lower than under inelastic prices.

\subsection{Additional Proofs}

\begin{proof}[Proof of Proposition~\ref{Prop:y}]
    The proof is similar to the one of Proposition~1 by \cite{V2GRO}, which analyzes the state-of-charge of an electric vehicle battery providing frequency regulation. The difference is that we do not model any electricity consumption for driving. By definition,
    \begin{align*}
    	y\big(x^b, x^r, \delta, y_0, t\big) = 
    	& \:
    	y_0 + \int_0^t \eta^+ \big[x^b + \delta(t') x^r \big]^+ - \frac{1}{\eta^-} \big[x^b + \delta(t') x^r \big]^- \, \mathrm{d}t'\\
    	= & \: y_0 + \int_0^t \min\left\{ \eta^+\big( x^b + \delta(t')x^r \big), \, \frac{1}{\eta^-} \big( x^b + \delta(t')x^r \big) \right\} \, \mathrm{d}t',
    \end{align*}
    where the second equality holds because $0 < \eta^+ \leq \nicefrac{1}{\eta^-}$. As $\eta^+ > 0$, $\eta^- >0$, and $x^r \geq 0$, both $\eta^+(x^b + \delta(t')x^r)$ and $\frac{1}{\eta^-}(x^b + \delta(t')x^r)$ are nondecreasing in~$\delta(t')$ and strictly increasing in~$x^b$. The minimum of two (nondecreasing/strictly increasing) affine functions is a concave (nondecreasing/strictly increasing) function~\citep[p.~73]{SB04}. The function~$y$ is thus concave strictly increasing in~$x^b$, concave in~$x^r$, concave nondecreasing in~$\delta$, and affine nondecreasing in~$y_0$.
    \end{proof}
    
\begin{proof}[Proof of Lemma~\ref{lem:D}]
    The claim follows immediately from the definition of~$\D$.
\end{proof}

\begin{proof}[Proof of Lemma~\eqref{eq:D_down}]
We first note that~$\D$ can be replaced with~$\D^+$ on the left-hand side of equation~\eqref{eq:D_down} because~$y$ is nondecreasing in~$\delta$ by Proposition~\ref{Prop:y} and because~$\D$ is symmetric by Lemma~\ref{lem:D}.

By construction, $\D^+_\downarrow \subseteq \D^+$. The claim thus follows if, for every signal~$\delta \in \D^+$, we can construct a rearranged signal~$\delta_\downarrow \in \D^+_\downarrow$ such that 
\begin{equation}\label{eq:soc_ordering}
    y(x^b, x^r, \delta, y_0, t) \leq 
    y(x^b, x^r, \delta_\downarrow, y_0, t)~~\forall t \in \T.
\end{equation}
We fix an arbitrary $\delta \in \D^+$ and define its rearrangement 
$ \delta_\downarrow(t) = \sup \{\rho \in [0,1] \, : \, \int_\T 1_{\delta(t') \geq \rho} \, \mathrm{d}t' \geq t\}$
involving the indicator function $1_{\delta(t') \geq \rho} = 1$ if $\delta(t') \geq \rho$, and $= 0$ otherwise. Note that the maximization problem is feasible because $\delta \in [0,1]$. We will now show that $\delta_\downarrow \in \D^+_\downarrow$.

To this end, we first note that, by definition, $\delta_\downarrow(t) \in [0,1]$ for all $t \in \T$. The signal~$\delta_\downarrow$ is also nonincreasing because the feasible set of the maximization problem becomes smaller as~$t$ grows. In addition, $\delta_\downarrow$ is left-continuous. To see this, consider the random variable $Y = \delta(t')$, where $t'$ follows the uniform distribution on~$T$. Then, $ \delta_\downarrow(t) = \sup \{\rho \in [0,1] \, : \, \mathbb{P}[Y \geq \rho] \geq \frac{t}{T}\}$, which is left-continuous because the function~$\mathbb{P}[Y \geq \rho]$ is left-continuous and nonincreasing.

Since~$\delta$ is Riemann integrable, it is also Lebesgue integrable. By the definition of the Lebesgue integral, we have 
\begin{equation*}
    \int_\T 1_{\delta(t) \in \set{B}} ~ \mathrm{d}t = \int_\T 1_{\delta_\downarrow(t) \in \set{B}} ~ \mathrm{d}t 
\end{equation*}
for any Borel set $\set{B} \subseteq [0, 1]$, which implies that 
\begin{equation*}
    \int_\T \delta_\downarrow(t) \, \mathrm{d}t = \int_\T \delta(t) \, \mathrm{d}t
    \leq \gamma,
\end{equation*}
where the inequality holds because $\delta \in \D^+$. In summary, we have thus shown that $\delta_\downarrow \in \D^+_\downarrow$.

We will now show the inequality in~\eqref{eq:soc_ordering}. The state-of-charge function can be written as
\begin{align*}
    y(x^b, x^r, \delta, y_0, t)
    = & y_0 + \int_0^t \min\left\{
    \eta^+ \left(x^b + \delta(t') x^r\right), \frac{1}{\eta^-} \left(x^b + \delta(t') x^r\right) \right\} \, \mathrm{d}t' \\
    = & y_0 + \int_0^T \min\left\{
    \eta^+ \left(x^b + \delta(t') x^r\right), \frac{1}{\eta^-} \left(x^b + \delta(t') x^r\right) \right\} \cdot 1_{t' \leq t} ~ \mathrm{d}t'
    \\
    \leq & y_0 + \int_0^T \min\left\{
    \eta^+ \left(x^b + \delta_\downarrow(t') x^r\right), \frac{1}{\eta^-} \left(x^b + \delta_\downarrow(t') x^r\right) \right\} \cdot 1_{t' \leq t} ~ \mathrm{d}t' \\
    = & y(x^b, x^r, \delta_\downarrow, y_0, t),
\end{align*}
where the first equality follows from the definition of~$y$ and from $0 < \eta^+ \leq \nicefrac{1}{\eta^-}$. The second equality follows from the definition of the indicator function. The inequality is a variant of the Hardy-Littlewood rearrangement inequality \citep[Theorem 378]{GH88}, which applies because $1_{t' \leq t}$ is nonincreasing and because the integrand $\dot y(\delta(t')) = \min\{
    \eta^+ (x^b + \delta(t') x^r), \frac{1}{\eta^-} (x^b + \delta(t') x^r) \}$
is nondecreasing in~$\delta(t')$ and thus admits the nonincreasing rearrangement $\dot y(\delta_\downarrow(t'))$.
\end{proof}

\begin{proof}[Proof of Lemma~\ref{lem:sym}]
	We first prove that the symmetry of~$\mathbb{P}_\xi$ implies that $F(z) + F(-z) = 1 + \mathbb{P}_\xi[\{z\}]$ for all $z \in \R$. To see this, note that
	\begin{equation*}\label{eq:Fsym}
		F(z) + F(-z) 
		= \mathbb{P}_\xi\left((-\infty, z]\right) + \mathbb{P}_\xi\left((-\infty,-z]\right)
		= \mathbb{P}_\xi\left((-\infty, z]\right) + \mathbb{P}_\xi\left([z,\infty)\right)
		= 1 + \mathbb{P}_\xi\left[\{z\}\right].
	\end{equation*}
	Thus, $h(z) = F(z) - 1/2(1 + \mathbb{P}_\xi[\{z\}]) = 1/2(1 + \mathbb{P}_\xi[\{z\}]) - F(-z) = -(F(-z) - 1/2(1 + \mathbb{P}_\xi[\{-z\}])) = -h(-z)$ is an odd function, which implies that
	\begin{align*}
		\varphi(z) & = \int_{-\infty}^{z} F(z') \, \mathrm{d} z'
		= \int_{-\infty}^{-z} F(z') \, \mathrm{d} z' + \int_{-z}^{z} F(z') - \frac{1}{2}(1 + \mathbb{P}_\xi[\{z'\}]) + \frac{1}{2}(1 + \mathbb{P}_\xi[\{z'\}]) \, \mathrm{d} \xi \\
		& = \varphi(-z) + \int_{-z}^{z} \frac{1}{2}(1 + \mathbb{P}_\xi[\{z'\}]) \, \mathrm{d} \xi
		= \varphi(-z) + z.
	\end{align*}
    The last equality holds as $\int_{-z}^{z} \mathbb{P}_\xi[\{z'\}]) \, \mathrm{d} z' = 0$. Indeed, the function $\mathbb{P}_\xi[\{z'\}])$ is nonzero only for countably many values of~$z'$ and is thus almost surely zero with respect to the Lebesgue~measure.
\end{proof}
\end{document}